\pdfoutput=1
\documentclass[runningheads]{llncs}
\usepackage{amssymb}
\usepackage{amsmath}
\usepackage{booktabs} % For pretty tables
\usepackage{caption} % For caption spacing
\usepackage{subcaption} % For sub-figures
\usepackage{graphicx}
\usepackage{pgfplots}
\usepackage[all]{nowidow}
\usepackage[utf8]{inputenc}
\usepackage{tikz}
\usetikzlibrary{er,positioning,bayesnet}
\usepackage{multicol}
\usepackage{algpseudocode,algorithm,algorithmicx}
\usepackage{multirow}
\usepackage[inline]{enumitem} % Horizontal lists
\usepackage{hyperref}

\usepackage{fourier} %Para el simbolo de peligro (aunque me cambio todo el estilo del papar, no me disgusta)
\usepackage{ upgreek }

\definecolor{blue}{HTML}{1F77B4}
\definecolor{orange}{HTML}{FF7F0E}
\definecolor{green}{HTML}{2CA02C}

\pgfplotsset{compat=1.14}

\allowdisplaybreaks

\begin{document}
\title{ Decision-dependent Wasserstein Distributionally Robust Optimization}
%Otra opción de titulo es: Distributionally Robust Optimization with ambiguity set as a Wasserstein ball with decision-dependent radius and center.

%
\titlerunning{Decision-dependent Distributionally Robust Optimization}
%Otra opcion de encabezado es: DRO using a Decision-Dependent Wasserstein Ball Ambiguity Set.
\author{Diego Fonseca \and
Mauricio Junca }
\institute{Department of Mathematics, Universidad de los Andes, Bogot\'a, Colombia
\email{\{df.fonseca,mj.junca20\}@uniandes.edu.co}}
\maketitle              % typeset the header of the contribution
\begin{abstract}
This work presents a new Distributionally Robust Optimization approach, using $p$-Wasserstein metrics, to analyze a stochastic program in a general context. The ambiguity set in this approach depends on the decision variable and is represented as a ball where both the center and the radius depend on the decision variable. We show that, under Lipschitz's assumptions for the objective function, our approach can be reformulated as a finite-dimensional optimization problem, which is sometimes convex. In addition, we numerically compare our proposed approach with the standard formulation of distributionally robust optimization, which typically does not use ambiguity sets dependent on the decision variable, in the context of portfolio optimization.  

\end{abstract}
\section{Introduction}

In this work, we consider general stochastic programs given by the formulation
\begin{equation}\label{StochsticProgWithExpectConst}
J={\displaystyle \min_{x\in\mathcal{X}} } \: \mathbb{E}_{\xi\sim\mathbb{P}}\left[F(x,\xi)\right],
\end{equation}
where $F$ is a function such that  $F:\mathbb{R}^{m}\times \mathbb{R}^{n}\rightarrow\mathbb{R}$, $\xi\in\mathbb{R}^{n}$ is a random vector with (unknown) probability distribution $\mathbb{P}$ supported in $\Xi\subseteq\mathbb{R}^{n}$, and $\mathcal{X}\subseteq \mathbb{R}^{m}$ is a set of constraints on the decision vectors. This problem is prevalent in various contexts, including finance \cite{Li2017,Dentcheva2003}, operations research \cite{Miller1965}, and machine learning \cite{Rigollet2011,Mu2017}. The significance of the problem (\ref{StochsticProgWithExpectConst}) stems from the fact that several functions that quantify risk or loss can be described as expected values. For instance, in a finance context, the Conditional Value-at-Risk (CVaR) is a function that can be represented as an expected value \cite{Rockafellar2000}. As such, it is feasible to reframe a wide range of problems that feature a CVaR as an objective function into the form of an optimization problem (\ref{StochsticProgWithExpectConst}).

The most widely used approach to solve (\ref{StochsticProgWithExpectConst})  is Sample Average Approximation (SAA), where samples of $\xi$ are used to replace the expected values with their sample mean. In recent years, there has been a focus on using stochastic gradient descent methods to solve the problem (\ref{StochsticProgWithExpectConst}). Examples of this include \cite{Akhtar2021,Lan2020,Xiao2019}. However, these methods can be sensitive to alterations in sample quality, and the out-of-sample performance may be poor, especially when the sample size is small.

We propose a data-driven approach for addressing (\ref{StochsticProgWithExpectConst}) using Distributionally Robust Optimization (DRO). In that sense, the DRO approach for the problem (\ref{StochsticProgWithExpectConst}) is formulated as
\begin{equation}\label{MarkovizRobustGeneral}
J_{\mathcal{D}}:=\min_{x\in\mathbb{X}}\sup_{\mathbb{Q}\in\mathcal{D}}\mathbb{E}_{\xi\sim\mathbb{Q}}[F(x,\xi)],
\end{equation}
where $\mathcal{D}$ is a set of probability distributions, which is known as \textit{ambiguity set}. Note that $J^{*}\leq J_{\mathcal{D}}$ if $\mathbb{P}\in \mathcal{D}$.  The set $\mathcal{D}$ plays a crucial role in the tractability of the problem under consideration. There have been several proposals in the literature on how to define $\mathcal{D}$. For instance, in references \cite{Lagoa2002} and \cite{Shapiro2006}, it is defined as a set of distributions supported at a single point. In contrast, references \cite{Delage2010,Popescu2007,Scarf1958,Shapiro2002} define $\mathcal{D}$ as the set of distributions satisfying certain moment restrictions or belonging to a particular parametric family of distributions.

Another option is to endow the set of probability distributions with a notion of distance, and define $\mathcal{D}$ as a ball in this metric. This ball is often centered on an empirical distribution $\widehat{\mathbb{P}}_{N}$, computed from a sample $\widehat{\xi}_{1},\ldots,\widehat{\xi}_{N}$ of the random vector $\xi$, with the radius chosen such that $\mathbb{P}$ belongs to the ball with high probability or such that the out-of-sample performance of the optimal solution is satisfactory. The choice of the distance metric influences the tractability of the resulting DRO. Commonly used metrics include Burg's entropy \cite{Wang2015}, Kullback-Leibler divergence \cite{Jiang2015}, and Total Variation distance \cite{Sun2015}. In this work, we adopt the Wasserstein distance and define $\mathcal{D}$ as a ball in this metric centered on the empirical distribution and with a properly chosen radius. Note that if the radius is set to 0 in this approach, we recover the SAA strategy.

\begin{definition}[Wasserstein distance]\label{Def:MetricaWasserstein}
The \textit{Wasserstein distance} $W_{p}(\mu,\nu)$ between $\mu,\nu\in\mathcal{P}_{p}(\Xi)$ is defined by
\begin{equation*}
\resizebox{.98 \textwidth}{!}{${\displaystyle W_{p}(\mu,\nu):=\left(\inf_{\Pi\in\mathcal{P}(\Xi\times\Xi)}\left\{\int_{\Xi\times\Xi}\mathbf{d}^{p}(\xi,\zeta)\Pi(d\xi,d\zeta)\: :\: \Pi(\cdot \times\Xi)=\mu(\cdot),\: \Pi(\Xi\times\cdot)=\nu(\cdot)\right\}\right)^{1/p}
}$}
\end{equation*}
where 
$$\mathcal{P}_{p}(\Xi):=\left\{\mu\in\mathcal{P}(\Xi)\: :\: \int_{\Xi}\mathbf{d}^{p}(\xi,\zeta_{0})\mu(d\xi) < \infty\ \mbox{for some }\zeta_{0}\in\Xi\right\}$$
and $d$ is a metric in  $\Xi$.
\end{definition}
$W_{p}$ defines a metric in $\mathcal{P}_{p}(\Xi)$ for  $p\in[1,\infty)$, hence, the  ball with respect to some  $p$-Wasserstein distance with radius $\varepsilon>0$ and center $\mu\in \mathcal{P}(\Xi)$ is given by
\begin{equation}\label{BolaRespectoP}
\mathcal{B}_{\varepsilon}\left(\mu\right):=\left\{\nu\in \mathcal{P}(\Xi) \: \left|\: W_{p}(\mu,\nu)\leq\varepsilon \right.\right\}.
\end{equation}

The Wasserstein distance, also referred to as Earth's moving distance in computer science, the Monge-Kantorovich-Rubinstein distance in physics, and the Optimal Transport distance in optimization, was first defined in \cite{Vasershtein1969}. Although it arose in various fields of science almost simultaneously, it is known by different names depending on the context.

There are numerous theoretical and practical reasons that make the Wasserstein distance particularly appealing, as highlighted in \cite{Villani2008}. One of its key advantages is its dual representation, which enables a more tractable equivalent formulation of (\ref{MarkovizRobustGeneral}). Specifically, using $p$-Wasserstein distances we obtain the following problem
\begin{equation} \label{DROWGeneral}
\widehat{J}^{\mathrm{S}}_{N,p,q}(\varepsilon):=\min_{x\in\mathbb{X}}\sup_{\mathbb{Q}\in\mathcal{B}_{\varepsilon}\left(\widehat{\mathbb{P}}_{N}\right) }\mathbb{E}_{\xi\sim\mathbb{Q}}[F(x,\xi)],
\end{equation}
where $\widehat{\mathbb{P}}_{N}$ is the empirical measure generated by sample $\widehat{\xi}_{1},\ldots,\widehat{\xi}_{N}$ of $\xi$, and the ball $\mathcal{B}_{\varepsilon}\left(\widehat{\mathbb{P}}_{N}\right)$ is defined with respect to the $p$-Wasserstein distance in $\mathbb{R}^{n}$ where the cost function used is $\mathbf{d}=\left\|\cdot\right\|_{q}$ (see definition \ref{Def:MetricaWasserstein}). In this work, we employ the notation "S" to denote "standard." This convention is adopted to highlight that the original formulation using Wasserstein's distance for distributionally robust optimization problems, referred to as the "standard" formulation, was the first of its kind. Additionally, this notation allows differentiating the standard formulation from our subsequent proposal. Finally, the problem (\ref{DROWGeneral}) can be reformulated by the following theorem.

\begin{theorem}\label{Thm:ReformulacionDROWInterno}
Assume that $ F $  is upper semicontinuous with respect to $\xi$. Then the problem (\ref{DROWGeneral}) is equivalent to the optimization problem
\begin{equation}\label{Eqn:ReformulacionDROW}
\left\{
\begin{array}{lll}
{\displaystyle \inf_{x\in\mathbb{X},\lambda,s}} & {\displaystyle \lambda \varepsilon^{p} +\frac{1}{N}\sum_{i=1}^{N}s_{i}} &\\
\mbox{subject to} & {\displaystyle \sup_{\xi\in\Xi}\left(f(x,\xi)-\lambda d^{p}(\xi,\widehat{\xi}_{i}) \right) \leq s_{i}  } & \forall i=1,\ldots,N, \\
&\lambda \geq 0.& 
\end{array}
\right.
\end{equation}
\end{theorem}

The previous theorem is formulated and proved in \cite{Blanchet2019}. However, the reformulation (\ref{Eqn:ReformulacionDROW}) has also been obtained under more restrictive assumptions in \cite{Kuhn2018} and \cite{Mehrotra2017}.

In many cases, depending on the form of the function $F$, the problem (\ref{Eqn:ReformulacionDROW}) may result in a semi-infinite optimization problem with a large number of variables, which can pose a significant challenge. This is because the supremum appearing in the constraints of (\ref{Eqn:ReformulacionDROW}) may not be solvable explicitly, and it is well-known that solving semi-infinite programs is computationally demanding. Thus, our aim is to obtain a reformulation of (\ref{MarkovizRobustGeneral}) that is not a semi-infinite problem. This motivates us to seek a different alternative ambiguity set that leads to a problem that is computationally more tractable.

Our proposal is based on the use of an ambiguity set as a ball in accordance with (\ref{BolaRespectoP}), but with a radius and center that depend on the decision variable. This is what is referred to as a decision-dependent ambiguity set, which will be further explained in Section \ref{Sec:MeanVariance}. The use of decision-dependent ambiguity sets using the Wasserstein metric has been explored in previous studies, such as \cite{Luo2020} and \cite{Noyan2021}. However, our proposed ambiguity set differs from the ones presented in these studies. In \cite{Luo2020}, only the radius of the ball depends on the decision, with the dependency being unspecified. In contrast, in \cite{Noyan2021}, the probability distributions within the ball depend on the decision, but the radius remains constant.

Finally, although the  primary focus of this work is to address the problem (\ref{StochsticProgWithExpectConst}),  a small portion of this work is dedicated to minimizing variance instead of minimizing expected value. Although a proposal similar to the one presented to address (\ref{StochsticProgWithExpectConst}) will be given for the variance case, we will only show its reformulation. Few studies have been conducted on solving stochastic problems with variance as the objective function. Most existing studies have focused on the mean-variance problem, which is a distinct problem. 
%The mean-variance problem was first introduced in portfolio optimization \cite{Markowitz1952}. Early attempts to solve the mean-variance problem relied on estimates of the return vector and covariance matrix, which can be replicated in the problem (\ref{StochsticProgWithExpectConst}), but these performed poorly out-of-sample and were sensitive to variations in the sample that define the estimates \cite{Chopra1993}. To overcome this, various works have used sets based on prior information about the returns or computationally tractable sets \cite{El-Ghaoui2003,Zymler2011,Natarajan2010,Lotf2017,Lotf2018,Won2020}. However, imposing unverifiable assumptions about the moments of the returns can also negatively impact the out-of-sample performance of these methods. Mean-variance models also emerge in inventory management to address the newsvendor problem \cite{Chen200}, with relevant research works including \cite{Choi2008,Rubio-Herrero2015,ZhangSethi2020}.

Concretely, in this work, we propose a new DRO with the Wasserstein metric to address stochastic programs from a sample of the random vector. In that sense, our contributions are the following:
\begin{enumerate}
\item[$\bullet$] We demonstrate that the problem resulting from our proposed distributional robust approach can be reformulated as an optimization problem with finite-dimensional variables.
\item[$\bullet$]  We prove that, under specific conditions, our proposed distributional robust approach is equivalent to the standard approach described in (\ref{DROWGeneral}), which assumes the ambiguity set is independent of the decision variables.
\item[$\bullet$]  In the context of mean-risk portfolio optimization, we show that our proposed approach tends to offer more computationally tractable optimization problems compared to those obtained through the standard approach outlined in (\ref{DROWGeneral}). Furthermore, our numerical simulations demonstrate that the decision-dependent ambiguity set approach generally performs better than the standard approach when all problem variables are considered in performance evaluation.
\item[$\bullet$]  We propose an approach similar to the one in (\ref{MarkovizRobust}) but tailored to address a stochastic optimization problem where the objective function is the variance instead of the expected value. We prove that the resulting problem can be reformulated as an optimization problem with finite-dimensional variables.
\end{enumerate}

The organization of this paper is as follows. In Section \ref{Sec:MeanVariance}, using Wasserstein distance, we describe our distributionally robust optimization model. In Section \ref{Sec:ProblemReformulation}, we also derive tractable reformulations for the optimization problem generated for our approach. Finally, in Section 4, we present some numerical results in the context of mean-risk portfolio optimization with risk measure as the Conditional Value at Risk CVaR.  
%Finally, we comment on future work in Section \ref{Sec:Conclusions}.

\subsubsection*{Notation:} For  $q\in[1,\infty)\cap \mathbb{N}$, the $q$-norm in $\mathbb{R}^{k}$ is is noted as $\left\|\cdot\right\|_{q}$. For $N\in \mathbb{N}$, we let $[N] := \{1, 2,\ldots,N\}$. Additionally, the conjugate of the function $f:\mathbb{R}^{n}\rightarrow \mathbb{R}$ is $f^{*}(\xi):=\sup_{\zeta\in\mathbb{R}^{n}}(\langle \zeta,\xi\rangle-f(\zeta))$, and its $q$-Lipschitz norm $\left\|f\right\|_{\mathrm{Lip,q}}:=\sup_{x\neq y}(f(x)-f(y))/\left\|x-y\right\|_{q}$. Finally, we use $\widetilde{O}_{\mathbb{P}}$ for the big $O$ in probability notation suppressing the logarithmic dependence.

%***********************************************************************
\section{New DRO formulation} \label{Sec:MeanVariance}

As stated before, in (\ref{StochsticProgWithExpectConst}),  the distribution  $\mathbb{P}$ is unknown, and we assume we have access to realizations of the random vector $\xi$. That is, let  $\widehat{\xi}_{1},\ldots,\widehat{\xi}_{N}$ be a sample of $\xi$ which allows to estimate  $\mathbb{P}$ by means of the empirical distribution $\widehat{\mathbb{P}}_{N}:=\frac{1}{N}\sum_{i=1}^{N}\delta_{\widehat{\xi}_{i}}$. Note that if we replace $\mathbb{P}$ with $\widehat{\mathbb{P}}_{N}$ in  (\ref{StochsticProgWithExpectConst}), we obtain the SAA strategy, with the drawbacks already mentioned. Another possible approach could be to consider $\mathcal{D}=\mathcal{B}_{\varepsilon}\left(\widehat{\mathbb{P}}_{N}\right)$ as ambiguity set in (\ref{MarkovizRobustGeneral}). However, this strategy has disadvantages. Note that if we use Theorem \ref{Thm:ReformulacionDROWInterno} to reformulate the objective function in (\ref{MarkovizRobustGeneral}), for a general function $F$ we could obtain a complicated and non-tractable semi-infinite optimization problem. This motivates our proposal, in which we seek to choose an ambiguity set that allows us to obtain a reformulation of (\ref{MarkovizRobustGeneral}) that is tractable and with good performance.

To present our approach we impose the following assumption.

\begin{assumption}[Lipschitz]\label{AssumptionPrincipal}
We assume that $F$ is a $q$-Lipschitz function with respect to $\xi$. This is, for each $x$, there exists $\gamma_{x,F,q} >0$ such that $|F(x,\xi)-F(x,\zeta)|\leq \gamma_{x,F,q}\left\|\xi-\zeta\right\|_{q} $ for all $\xi,\zeta\in\Xi$. We denote $\gamma_{x,F,q}=\left\|F(x,\cdot)\right\|_{\mathrm{Lip,q}}$.
\end{assumption}

Our approach also uses an empirical distribution but this depends on $x$. First, we establish the following convention: For $x\in\mathbb{R}^{m}$, we define  $\zeta^{x,F}:=F\left( x,\xi \right)$, note that this is a  random variable. We called  $\mathbb{P}^{x,F}$   to the probability distribution of $\zeta^{x,F}$. Because it depends on  $\mathbb{P}$,  $\mathbb{P}^{x,F}$ is also unknown. Additionally, we define $\widehat{\zeta}^{x,F}_{i}:=F\left( x,\widehat{\xi}_{i}\right)$, so $\widehat{\zeta}^{x,F}_{1},\ldots,\widehat{\zeta}^{x,F}_{N}$ is a sample of $\zeta^{x,F}$. This allows us to define the empirical distribution of $ \zeta^{x, F}$, which is given by $\widehat{\mathbb{P}}^{x,F}_{N}:=\frac{1}{N}\sum_{i=1}^{N}\delta_{\widehat{\zeta}^{x,F}_{i}}$. Therefore, we consider the following optimization problem: For a given $\epsilon>0$

\begin{equation}\label{MarkovizRobust}
\widehat{J}^{\mathrm{A}}_{N,p,q}(\varepsilon):={\displaystyle \min_{x\in\mathcal{X}} }  {\displaystyle\sup_{\mathbb{Q}\in\mathcal{B}_{\varepsilon\gamma_{x,F,q}}\left(\widehat{\mathbb{P}}_{N}^{x,F}\right) } \mathbb{E}_{\zeta\sim\mathbb{Q}}[\zeta] }
\end{equation}
where $\mathcal{B}_{\varepsilon\gamma_{x,F,q}}\left(\widehat{\mathbb{P}}_{N}^{x,F}\right)$ is a ball centered at $\widehat{\mathbb{P}}_{N}^{x,F}$ with radius  $\varepsilon \gamma_{x,F,q}$. This ball is defined with respect to the $p$-Wasserstein distance in $\mathbb{R}$ where the cost function used is $\mathbf{d}=\left|\cdot\right|$. We adopt the letter "A" to refer to the term "alternative", with the purpose of distinguishing our proposal from the formulation presented in (\ref{DROWGeneral}).

The following result gives the reason why our ambiguity set has radius $\varepsilon \gamma_{x,F,q}$, and also shows a relationship between $\mathcal{B}_{\varepsilon\gamma_{x,F,q}}\left(\widehat{\mathbb{P}}_{N}^{x,F}\right)$ and $\mathcal{B}_{\varepsilon}(\widehat{\mathbb{P}}_{N})$.

\begin{lemma} \label{Lemma:JustificacionBola}
$ W_{p}(\widehat{\mathbb{P}}_{N}^{x,F},\mathbb{P}^{x,F}) \leq \gamma_{x,F,q}  W_{p}\left( \widehat{\mathbb{P}}_{N} ,\mathbb{P}\right)$ for $p,q\geq1$ where $ W_{p}\left( \widehat{\mathbb{P}}_{N} ,\mathbb{P}\right)$ is considered with cost function $\mathbf{d}=\|\cdot\|_{q}$ and $W_{p}(\widehat{\mathbb{P}}_{N}^{x,F},\mathbb{P}^{x,F}) $ with cost function $\mathbf{d}=|\cdot|$.
\end{lemma}

Note that if $\varepsilon>0$ is such that $\mathbb{P}\in\mathcal{B}_{\varepsilon}(\widehat{\mathbb{P}}_{N})$, where the ball is taken with respect to the $p$-Wasserstein metric  with cost function $\mathbf{d}=\|\cdot\|_{q}$ and distributions supported in a subset of  $\mathbb{R}^{n}$, then $\mathbb{P}^{x,F}\in\mathcal{B}_{\varepsilon \gamma_{x,F}} (\widehat{\mathbb{P}}_{N}^{x,F})$, where these balls are taken with respect to the $p$-Wasserstein metric with cost function $\mathbf{d}=|\cdot|$ and  distributions supported in a subset of $\mathbb{R}$. The proof of Lemma \ref{Lemma:JustificacionBola} is addressed in \ref{Apendice:PruebasLema}. 

As a direct result of Theorem \ref{Thm:ReformulacionDROWInterno}, we can deduce the following lemma, which enables us to generate an initial reformulation of (\ref{MarkovizRobust}).
\begin{lemma} \label{Lemma:ReformulationInicial} For each $p,q\geq 1$, $\varepsilon\geq 0$ and $x\in\mathcal{X}$,
    \begin{equation}\label{eqn:Lemma:ReformulacionPreliminar}
    \sup_{\mathbb{Q}\in\mathcal{B}_{\varepsilon\gamma_{x,F,q}}\left(\widehat{\mathbb{P}}_{N}^{x,F}\right) } \mathbb{E}_{\zeta\sim\mathbb{Q}}[\zeta] =\inf_{\lambda\geq 0}\left(\lambda \varepsilon^{p}\gamma_{x,F,q}^{p}+\frac{1}{N}\sum_{i=1}^{N}\sup_{\zeta\in F(x,\Xi)}\left(\zeta -\lambda\left|\zeta -F(x,\widehat{\xi}_{i})\right|^{p} \right)\right).
    \end{equation}
\end{lemma}
It should be noted that Lemma \ref{Lemma:ReformulationInicial} can be used to derive a reformulation of equation (\ref{MarkovizRobust}) by minimizing expression (\ref{eqn:Lemma:ReformulacionPreliminar}) with respect to the variable $x$. The reformulation of (\ref{MarkovizRobust}) resulting from this lemma is a preliminary one, which has the potential to be further enhanced regarding its computational tractability. Despite its preliminary nature, this lemma plays a crucial role in the following proposition, which establishes a relationship between (\ref{DROWGeneral}) and (\ref{MarkovizRobust}). The corollary of this proposition will further elucidate this relationship.

\begin{proposition} For each $p,q\geq 1$, $\varepsilon\geq 0$ and $x\in\mathcal{X}$,
$$ \sup_{\mathbb{Q}\in\mathcal{B}_{\varepsilon}\left(\widehat{\mathbb{P}}_{N}\right) }\mathbb{E}_{\xi\sim\mathbb{Q}}[F(x,\xi)] \leq \sup_{\mathbb{Q}\in\mathcal{B}_{\varepsilon\gamma_{x,F,q}}\left(\widehat{\mathbb{P}}_{N}^{x,F}\right) } \mathbb{E}_{\zeta\sim\mathbb{Q}}[\zeta].$$
\end{proposition}
\proof For each $x\in\mathcal{X}$, the proof is presented in the following lines:
\begin{align}
    \sup_{\mathbb{Q}\in\mathcal{B}_{\varepsilon}\left(\widehat{\mathbb{P}}_{N}\right) }\mathbb{E}_{\mathbb{Q}}[F(x,\xi)] &= \inf_{\beta\geq 0}\left(\beta \varepsilon^{p} +\frac{1}{N}\sum_{i=1}^{N}\sup_{\xi\in\Xi}\left(F(x,\xi)-\beta \left\|\xi-\widehat{\xi}_{i}\right\|_{q}^{p}  \right)\right) \label{PrrofInequaLine1}\\
    &\leq  \inf_{\beta\geq 0}\left(\beta \varepsilon^{p} +\frac{1}{N}\sum_{i=1}^{N}\sup_{\xi\in\Xi}\left(F(x,\xi)-\frac{\beta}{\gamma_{x,F,q}^{p}} \left|F(x,\xi)-F(x,\widehat{\xi}_{i})\right|^{p}  \right)\right) \label{PrrofInequaLine2}\\ 
    &=   \inf_{\lambda\geq 0}\left(\lambda \varepsilon^{p}\gamma_{x,F,q}^{p} +\frac{1}{N}\sum_{i=1}^{N}\sup_{\xi\in\Xi}\left(F(x,\xi)-\lambda \left|F(x,\xi)-F(x,\widehat{\xi}_{i})\right|^{p}  \right)\right) \label{PrrofInequaLine3}\\
    &=\inf_{\lambda\geq 0}\left(\lambda \varepsilon^{p}\gamma_{x,F,q}^{p}+\frac{1}{N}\sum_{i=1}^{N}\sup_{\zeta\in F(x,\Xi)}\left(\zeta -\lambda\left|\zeta -F(x,\widehat{\xi}_{i})\right|^{p} \right)\right) \label{PrrofInequaLine4}\\
    &= \sup_{\mathbb{Q}\in\mathcal{B}_{\varepsilon\gamma_{x,F,q}}\left(\widehat{\mathbb{P}}_{N}^{x,F}\right) } \mathbb{E}_{\zeta\sim\mathbb{Q}}[\zeta]. \label{PrrofInequaLine5}
\end{align}
Equality (\ref{PrrofInequaLine1}) is a direct consequence of Theorem \ref{Thm:ReformulacionDROWInterno}, while inequality (\ref{PrrofInequaLine2}) is established based on Assumption \ref{AssumptionPrincipal}. Equality (\ref{PrrofInequaLine3}) is derived from the change of variable $\lambda= \frac{\beta}{\gamma_{x,F,q}^{p}}$, and equality \ref{PrrofInequaLine4}) results from the change of variable $\zeta=F(x,\xi)$. Finally, equality (\ref{PrrofInequaLine5}) is a direct consequence of Lemma \ref{Lemma:ReformulationInicial}.
\qed
\endproof

\begin{corollary}
    $\widehat{J}^{\mathrm{S}}_{N,p,q}(\varepsilon)\leq \widehat{J}^{\mathrm{A}}_{N,p,q}(\varepsilon) $ for $p,q\geq1$ and $\varepsilon\geq 0$.
\end{corollary}

From now on, the optimal solutions of (\ref{DROWGeneral}) and (\ref{MarkovizRobust}) will be denoted as $\widehat{x}_{N,p,q}^{\mathrm{S}}(\varepsilon)$ and $\widehat{x}_{N,p,q}^{\mathrm{A}}(\varepsilon)$  respectively.

%*************************************************
\section{Problem reformulation} \label{Sec:ProblemReformulation}

In this section, the objective is to reformulate (\ref{MarkovizRobust}) as an optimization problem with finite-dimensional variables. This reformulation depends on the image of the support of $\xi$ under the function $F(x,\cdot)$ for each $x\in\mathcal{X}$.

\begin{theorem} \label{Prop:Reformul1MarkovizDROWReformLargaLargaSupportAcot} \hfill
\begin{enumerate}
    \item[(a)] If $p=1$ and $F(x,\Xi)$ is an interval for each $x\in\mathcal{X}$, then   the optimization problem (\ref{MarkovizRobust})  is equivalent to the following optimization problem
\begin{equation}\label{MarkovizDROWReformLargaSupportAcot}
    \widehat{J}_{N,p,q}^{\mathrm{A}}(\varepsilon)= \underset{ x\in\mathcal{X}  }{\mathrm{minimize}} \min\left\{\frac{1}{N}{\displaystyle\sum_{i=1}^{N}}F\left(x,\widehat{\xi}_{i}\right)+\varepsilon\gamma_{x,F,q},\sup_{\xi\in\Xi}F(x,\xi)\right\}.
\end{equation}

\item[(b)] If $p\geq 1$, $F(x,\Xi)$ is an interval, and $\sup_{\xi\in\Xi}F(x,\xi)=\infty$ for each  $x\in\mathcal{X}$, then the optimization problem (\ref{MarkovizRobust})  is equivalent to the following optimization problem
\begin{equation}\label{MarkovizDROWReformLargaSupportNoAcot}
     \widehat{J}_{N,p,q}^{\mathrm{A}}(\varepsilon)=\underset{ x\in\mathcal{X}  }{\mathrm{minimize}}   \frac{1}{N}{\displaystyle\sum_{i=1}^{N}}F\left(x,\widehat{\xi}_{i}\right)+\varepsilon\gamma_{x,F,q}\left(\frac{1}{p}+\frac{p-1}{p^{1/(p-1)}}\right).
\end{equation}

\end{enumerate}
\end{theorem}

In principle, case (a) is considered less restrictive compared to case (b) since it demands fewer conditions on the image of support $\Xi$ under $F(x,\cdot)$. Specifically, in case (a), for a given $x \in \mathcal{X}$, $F(x,\Xi)$ may be a bounded or unbounded interval, whereas in case (b), $F(x,\Xi)$ is required to be unbounded. However, it is worth noting that there are several crucial functions that satisfy it. For instance, for $\Xi=\mathbb{R}^{m}$ the function $F(x,\xi):=\left\langle x,\xi\right\rangle$ satisfies this condition and is particularly significant in the context of portfolio optimization. Additionally, any function that can be expressed as the maximum of a finite number of affine linear functions with respect to $\xi$ also satisfies the conditions of the case (b) for $\Xi=\mathbb{R}^{m}$.

\remark{(\textit{Relationship with the regularizing approach.})} \label{RemarkRegularizacion}
Under certain conditions on the support $\Xi$, the optimization problems (\ref{MarkovizDROWReformLargaSupportAcot}) and (\ref{MarkovizDROWReformLargaSupportNoAcot}) can be expressed as the regularized sampling problem, and its use in the estimation of (\ref{DROWGeneral}) has been previously explored. For certain objective functions $F$, (\ref{DROWGeneral}) is exactly a regularized sampling problem, as demonstrated in \cite{Shafieezadeh-AbadehKuhn2019,Blanchet2017}. However, the objective functions considered in these studies focus on regression and classification problems in the context of supervised learning, which requires a special cost function for the 1-Wasserstein distance that they adopt in their results. In contrast, a more general context is considered in \cite{GaoChenKleywegt2017}, where a regularized version of the problem is proposed for $p\in [1,\infty)\cap \mathbb{N}$, and an asymptotic bound of the distance between the optimal values of (\ref{DROWGeneral}) and the regularized sampling problem is established with respect to the sample size.

The main result of \cite{GaoChenKleywegt2017} is to show that
$$
{\displaystyle\sup_{\mathbb{Q}\in\mathcal{B}_{\varepsilon}\left(\widehat{\mathbb{P}}_{N} \right) } \mathbb{E}_{\zeta\sim\mathbb{Q}}\left[F(x,\zeta)\right] }\leq \mathbb{E}_{\widehat{\mathbb{P}}_{N}}[F(x,\xi)]+\varepsilon_{N}V(F(x,\cdot))+\widetilde{O}_{\mathbb{P}}\left(\frac{1}{N}\right) \quad \mbox{for each }x\in\mathcal{X}
$$
where $\varepsilon_{N}=\widetilde{O}_{\mathbb{P}}\left(\frac{1}{\sqrt{N}}\right)$, and $V$ is a function referred to as the \textit{variation of the loss}, which serves as the regularizing term. Specifically, for $p=1$, $V$ is $\gamma_{x,F,q}$, and for $p>1$, $V$ is defined as $\mathbb{E}_{\widehat{\mathbb{P}}_{N}}\left[\left|\nabla_{\xi}F(x,\xi) \right|_{\bar{p}}\right]^{\frac{1}{\bar{p}}}$, where $\bar{p}$ is the Hölder conjugate of $p$. The approach proposed in \cite{GaoChenKleywegt2017} and our approach differ completely for $p>1$, as guaranteed by Theorem \ref{Prop:Reformul1MarkovizDROWReformLargaLargaSupportAcot} (b).

\endremark

One of the concerns is to know when the new approach proposed in this paper is equal to the standard DRO in which the ball does not depend on the decision variable. This is clarified in the following result.

\begin{proposition} \label{Prop:IgualdadEnReformulaciones}
Suppose that $\Xi=\mathbb{R}^{d}$, $F(x,\xi)$ is a convex function respect to $\xi$,  and  $F(x,\Xi)$ is an interval with $\sup_{\xi\in\Xi}F(x,\xi)=\infty$ for each $x\in\mathcal{X}$. Then for each $x$ where $F(x,\cdot)$ is defined, 
$$ {\displaystyle\sup_{\mathbb{Q}\in\mathcal{B}_{\varepsilon\gamma_{x,F,q}}\left(\widehat{\mathbb{P}}_{N}^{x,F}\right) } \mathbb{E}_{\zeta\sim\mathbb{Q}}\left[\zeta\right] }=   {\displaystyle\sup_{\mathbb{Q}\in\mathcal{B}_{\varepsilon}\left(\widehat{\mathbb{P}}_{N} \right) } \mathbb{E}_{\xi\sim\mathbb{Q}}\left[F(x,\xi)\right] }
$$
if  any of the following conditions are satisfied:
\begin{enumerate}
    \item[(i)] $p=1$ and $q=2$.
    \item[(ii)] $p=2$, $q=2$, and $\sup_{z\in \Uptheta_{x,F} } \left( \frac{\|z\|^{2}_2}{4\lambda }+\langle z, \zeta \rangle-I_{x,F}^{*}(z) \right)= \frac{\gamma_{x,F,2}^{2} }{4\lambda}+I_{x,F}(\zeta) $ for all $x\in\mathcal{X}$, $\zeta\in \Xi$, and $\lambda >0$, where $I_{x,F}(\zeta):=F(x,\zeta)$, $I_{x,F}^{*}$ is the conjugate function of $I_{x,F}$, and  $\Uptheta_{x,F}:=\left\{z\in\mathbb{R}^{d}\: :\: I^{*}_{x,F}(z)<\infty \right\}$.
\end{enumerate}
\end{proposition}

Although condition (ii) appears to be restrictive, there are important functions that satisfy it. For example, as in Remark \ref{RemarkRegularizacion}, the function $F(x,\xi):=\langle x,\zeta \rangle$ satisfies (ii). Another aspect to highlight of this proposition is the type of support for which the equality is valid. This result shows that the support should be $\mathbb{R}^{d}$, which may be restrictive. However, empirical evidence suggests that for specific $F$ functions and condition (i), equality between the two approaches can be attained when supports are specified as convex polyhedra. 

\remark
The proof of the equivalence between (\ref{MarkovizRobust}) and (\ref{DROWGeneral}) for $p=1$ has been demonstrated in \cite{GaoChenKleywegt2017}. However, this proof is subject to different conditions than those outlined in Proposition \ref{Prop:IgualdadEnReformulaciones}.  Indeed, the result is formulated where the cost function is any norm, and the Lipschitz norm of $F(x,\cdot)$ must be reached at infinity, specifically, there must exist $\xi_{0}\in\Xi$ such that $\gamma_{x,F,q}=\limsup_{\|\xi-\xi_{0}\|_{q}\rightarrow \infty}\frac{F(x,\xi)-F(x,\xi_{0})}{\|\xi-\xi_{0}\|_{q}}$. 
\endremark

\subsection{Variance case}

In this part, we consider the scenario where the variance serves as an objective function in the stochastic optimization problem. To tackle this issue, we adopt the same method put forward for the problem (\ref{StochsticProgWithExpectConst}) and use Assumption \ref{AssumptionPrincipal}. Accordingly, the focus of this section is the following problem: 
$$J^{\mathrm{Var}}=\inf_{x\in\mathcal{X}}\mathrm{Var}_{\mathbb{P}}\left[F(x,\xi)\right].$$
The standard distributionally robust version using Wasserstein distances for this problem is of the form
\begin{equation} \label{VarianceRobust}
J_{N,p,q}^{\mathrm{Var},\mathrm{S}}(\varepsilon)=\inf_{x\in\mathcal{X}} {\displaystyle\sup_{\mathbb{Q}\in\mathcal{B}_{\varepsilon}\left(\widehat{\mathbb{P}}_{N}\right) }}\mathrm{Var}_{\xi\sim\mathbb{Q}}\left[F(x,\xi)\right].
\end{equation}
Our examination of current literature on the subject did not find works that endeavor to reformulate this problem in a general context. This dearth of exploration could potentially be attributed to the intricate nature that the inclusion of variance introduces in contrast to the expected value. Following the spirit of (\ref{MarkovizRobust}), the distributional robust version with dependent ambiguity set using Wasserstein distances for this problem is of the form
\begin{equation} \label{VarianceRobustDepend}
J_{N,p,q}^{\mathrm{Var},\mathrm{A}}(\varepsilon)=\inf_{x\in\mathcal{X}} {\displaystyle\sup_{\mathbb{Q}\in\mathcal{B}_{\varepsilon\gamma_{x,F,q}}\left(\widehat{\mathbb{P}}_{N}^{x,F}\right) }}\mathrm{Var}_{\zeta\sim\mathbb{P}}\left[\zeta\right].
\end{equation}
The following theorem states a reformulation of (\ref{VarianceRobustDepend}).

\begin{theorem}\label{Prop:Reformul1MarkovizDROWReformLargaLarga}
Assuming $p=q=2$, if $F(x,\Xi)=[0,\infty)$ or $F(x,\Xi)=\mathbb{R}$ for each $x\in\mathcal{X}$, then  the optimization problem (\ref{VarianceRobustDepend})  is equivalent to the following optimization problem with finite-dimensional variables
\begin{equation}\label{MarkovizDROWReformLarga}
     J_{N,p,q}^{\mathrm{Var},\mathrm{A}}(\varepsilon)=\underset{ x\in\mathcal{X}  }{\mathrm{minimize}}   \left( \sqrt{\frac{1}{N}{\displaystyle\sum_{i=1}^{N}}F\left(x,\widehat{\xi}_{i}\right) ^{2}-\frac{1}{N^{2}}\left({\displaystyle\sum_{i=1}^{N}}F\left( x,\widehat{\xi}_{i}\right)\right)^{2} }+\varepsilon\gamma_{x,F} \right)^{2}. 
     \end{equation}
\end{theorem}

%Given that reformulating (\ref{VarianceRobust}) is recognized as a challenging task, we can reasonably anticipate that reformulating (\ref{VarianceRobustDepend}) would also be arduous. Consequently, 

Theorem \ref{Prop:Reformul1MarkovizDROWReformLargaLarga} provides the reformulation exclusively for the scenario in which $p=q=2$. Extending this result to other values of $p$ and $q$ represents an unresolved challenge, as does identifying the conditions under which the standard approach is equivalent to the alternative approach. Notably, it is worth mentioning that the equality between the two approaches can be inferred from the outcomes presented in \cite{Blanchet2018} when $F(x,\xi)=\langle x,\xi\rangle$. 
%Nevertheless, for others functions, the question of whether the two approaches are equivalent remains unanswered.

%******************************************************************************************************

\section{ Mean-risk portfolio optimization}

In this section, we conduct numerical experiments in the context of mean-risk portfolio optimization. Our objective is  to illustrate that the approach proposed in this work may offer superior computational and performance advantages compared to the standard approach in certain cases. For this purpose, we replicate the framework proposed in Section 7 of \cite{Kuhn2018}. The problem addressed herein pertains to portfolio optimization, and its formulation is given by 
\begin{equation} \label{eqn:Mean-CVaROriginal}
    J=\inf_{x\in\mathcal{X}}\left(\mathbb{E}_{\xi\sim\mathbb{P}}[-\langle x, \xi\rangle]+\rho\mathrm{CVaR}_{\alpha,\xi\sim\mathbb{P}}(-\langle x,\xi\rangle)\right)
\end{equation}
where $\xi\in \Xi\subset\mathbb{R}^{m}$, and $\mathcal{X}:=\left\{x\in\mathbb{R}^{m}_{+}\: :\:\langle x, e \rangle=1\right\}$ with $e\in\mathbb{R}^{m}$ denoting the vector with all components equal to one. The constants $\rho\in\mathbb{R}$ and $\alpha\in(0,1]$ are determined by the investor and reflect the investor's risk aversion, and the parameter that makes $\mathrm{CVaR}_{\alpha,\xi\sim\mathbb{P}}$ the average value of the $\alpha$ percent of the most severe portfolio losses, respectively, under the probability distribution $\mathbb{P}$. The results presented in \cite{Rockafellar2000} enable us to express this problem in the form of the equation (\ref{StochsticProgWithExpectConst}), leading to the following equivalent formulation: 
\begin{equation} \label{eqn:Mean-CVaROriginalReformu}
J=\inf_{x\in\mathcal{X},\tau\in\mathbb{R}}\mathbb{E}_{\xi\sim\mathbb{P}}\left[\max_{k\in\left\{1,2\right\}}\left(a_{k}\langle x,\xi \rangle+b_{k}\tau\right)\right]
\end{equation}
where $a_{1}=-1$, $a_{2}=-1-\frac{\rho}{\alpha}$, $b_{1}=\rho$, and $b_{2}=\rho\left(1-\frac{1}{\alpha}\right)$.  The standard distributionally robust approach for this problem, referred to hereafter as WDROS and described in (\ref{DROWGeneral}), is the following:
\begin{equation} \label{CVaRWassertsStandard}
\widehat{J}^{\mathrm{S}}_{N,p,q}(\varepsilon) = \inf_{x\in\mathcal{X},\tau\in\mathbb{R}}\sup_{\mathbb{Q}\in\mathcal{B}_{\varepsilon}\left(\widehat{\mathbb{P}}_{N}\right) } \mathbb{E}_{\xi\sim\mathbb{Q}}\left[\max_{k\in\left\{1,2\right\}}\left(a_{k}\langle x,\xi \rangle+b_{k}\tau\right)\right].
\end{equation}
The optimal solution of this problem is denoted as $\widehat{x}_{N,p,q}^{\mathrm{S}}(\varepsilon)$ and $\widehat{\tau}_{N,p,q}^{\mathrm{S}}(\varepsilon)$. In the same direction, considering $F(x,\xi):=\max_{k\in\left\{1,2\right\}}\left(a_{k}\langle x,\xi \rangle+b_{k}\tau\right)$, the alternative distributively robust approach for this problem, referred to hereafter as WDROA, is given by:
\begin{equation} \label{CVaRWassertsAlternative}
    \widehat{J}^{\mathrm{A}}_{N,p,q}(\varepsilon)={\displaystyle \inf_{x\in\mathcal{X},\tau\in\mathbb{R}} }  {\displaystyle\sup_{\mathbb{Q}\in\mathcal{B}_{\varepsilon\gamma_{x,F,q}}\left(\widehat{\mathbb{P}}_{N}^{x,F}\right) } \mathbb{E}_{\zeta\sim\mathbb{Q}}[\zeta] }.
\end{equation}
The optimal solution of this problem is denoted as $\widehat{x}_{N,p,q}^{\mathrm{A}}(\varepsilon)$ and $\widehat{\tau}_{N,p,q}^{\mathrm{A}}(\varepsilon)$.

To facilitate comprehension of the results presented in this section, we must introduce several key concepts. The first is commonly referred to as out-of-sample performance. We will consider two types of out-of-sample performance in this study. The first type considers only the optimal decision associated with the variable $x$. That is, given a portfolio $\widehat{x}_{N}$ induced by the sample via some mechanism, for instance, $\widehat{x}_{N}$ can be $\widehat{x}_{N,p,q}^{\mathrm{S}}(\varepsilon)$ or $\widehat{x}_{N,p,q}^{\mathrm{A}}(\varepsilon)$, we define the out-of-sample performance of $\widehat{x}_{N}$ as 
\begin{equation} \label{eqn:DefOutOfSamplePerformSinTau}
J(\widehat{x}_{N}):=\mathbb{E}_{\xi\sim\mathbb{P}}[-\langle \widehat{x}_{N}, \xi\rangle]+\rho\mathrm{CVaR}_{\alpha,\xi\sim\mathbb{P}}(-\langle \widehat{x}_{N},\xi\rangle)=\inf_{\tau\in\mathbb{R}}\mathbb{E}_{\xi\sim\mathbb{P}}\left[\max_{k\in\left\{1,2\right\}}\left(a_{k}\langle \widehat{x}_{N},\xi \rangle+b_{k}\tau\right)\right].
\end{equation}
The objective is to identify $\widehat{x}_{N}$ decisions that minimize this function as much as possible. We refer to this type of out-of-sample performance as \textit{portfolio-focused out-of-sample performance}. 

We introduce the other type of out-of-sample performance, which considers both variables $x$ and $\tau$. Specifically, let $\widehat{x}_{N}$ and $\widehat{\tau}_{N}$ be the decisions induced by the sample through some mechanism. In this case, we define the out-of-sample performance of $\widehat{x}_{N}$ and $\widehat{\tau}_{N}$ as 
$$
J(\widehat{x}_{N},\widehat{\tau}_{N}):=\mathbb{E}_{\xi\sim\mathbb{P}}\left[\max_{k\in\left\{1,2\right\}}\left(a_{k}\langle \widehat{x}_{N},\xi \rangle+b_{k}\widehat{\tau}_{N}\right)\right].
$$
Again, the objective here is also to identify $\widehat{x}_{N}$ and $\widehat{\tau}_{N}$ decisions that minimize this function as much as possible. Hereafter, we refer to this type of out-of-sample performance as \textit{comprehensive out-of-sample performance}.

The two methods for evaluating out-of-sample performance are conceptually distinct. The first approach employs the $\tau$ corresponding to the true probability, which reflects the portfolio's performance in a setting that is unaffected by the sample. In contrast, the other method evaluates the out-of-sample performance by incorporating the $\tau$ generated by the sample, $\widehat{\tau}_{N}$. This method is more closely aligned to measure the degree to which $\widehat{x}_{N}$ and $\widehat{\tau}_{N}$ approximate the true solutions of  (\ref{eqn:Mean-CVaROriginalReformu}).

The last concept that we define is known as reliability which is linked to the concept of out-of-sample performance, so we have two types of reliability. The first one we will call \textit{portfolio-focused reliability}. To define this concept, we will focus momentarily on the standard case, the alternative case is analogous. Hence,  portfolio-focused reliability for the standard case is defined as $\mathbb{P}^{N}\left[J\left(\widehat{x}_{N,p,q}^{\mathrm{S}}(\varepsilon)\right) \leq \widehat{J}^{\mathrm{S}}_{N,p,q}(\varepsilon)\right]$. Analogously, the second reliability we will call \textit{comprehensive reliability}. Again, focusing on the standard case, this is defined as $\mathbb{P}^{N}\left[J\left(\widehat{x}_{N,p,q}^{\mathrm{S}}(\varepsilon),\widehat{\tau}_{N,p,q}^{\mathrm{S}}(\varepsilon)\right) \leq \widehat{J}^{\mathrm{S}}_{N,p,q}(\varepsilon)\right]$. This concept quantifies the probability that the optimal value of the distributionally robust approach exceeds the associated out-of-sample performance, which is something that should be available if the true probability distribution governing the random vector belongs to the Wasserstein ball associated with the approach under study. However, this is not a necessary condition for good out-of-sample performance. Finally, the portfolios generated by the two approaches studied will be compared with the $\infty$-norm in $\mathbb{R}^{m}$, in order to measure the differences component by component.

In light of the aforementioned considerations, we will examine two distinct scenarios. The first of these we shall refer to as the \textit{non-limited losses assets}, wherein we suppose that $\Xi=\mathbb{R}^{m}$. In this case, the results presented are the product of numerical experiments.  The second scenario we will investigate is referred to as the \textit{limited losses assets}, in which we assume that $\Xi=\left\{\xi\in\mathbb{R}^{m}\: :\:\xi\geq -1\right\}$. This last support reflects that no asset in the portfolio can lose more than 100\% of its value. For this last case, the results are theoretical.

\subsection{Non-limited-loss assets}

In this part, we assume that the random variable $\xi$ can be expressed as $\xi=\psi+\zeta_{i}$   where $\psi\sim \mathcal{N}(0,2\%)$ and $\zeta\sim\mathcal{N}(i\times 3\%, i\times 2.5\%)$ for each $i=1,\dots,N$, leading to the support of the random vector of returns given by $\Xi=\mathbb{R}^{m}$. We adopt the $p$-Wasserstein distance with $p=2$ and the cost function $\mathbf{d}=\|\cdot\|_{2}$. The standard distributionally robust approach that employs this Wasserstein distance,  described in (\ref{CVaRWassertsStandard}), can be reformulated as follows:
\begin{align}
    \widehat{J}^{\mathrm{S}}_{N,p,q}(\varepsilon)
    &= \left\{
\begin{array}{lll}
{\displaystyle \inf_{x\in\mathcal{X},\tau\in\mathbb{R},\lambda\geq 0,s\in\mathbb{R}^{N}}} & {\displaystyle \lambda \varepsilon^{2} +\frac{1}{N}\sum_{i=1}^{N}s_{i}} &\\
\mbox{subject to} & {\displaystyle \sup_{\xi\in\mathbb{R}^{m}}\left(\max_{k\in\left\{1,2\right\}}\left(a_{k}\langle x,\xi \rangle+b_{k}\tau\right)-\lambda \left\|\xi-\widehat{\xi}_{i} \right\|_{2}^{2} \right) \leq s_{i} , } & \forall i\in[N]. 
\end{array}
\right.  \label{eqn1:ReformulCVaRCasoP2} \\
&= \left\{
\begin{array}{lll}
{\displaystyle \inf_{x\in\mathcal{X},\tau\in\mathbb{R},\lambda\geq 0,s\in\mathbb{R}^{N}}} & {\displaystyle \lambda \varepsilon^{2} +\frac{1}{N}\sum_{i=1}^{N}s_{i}} &\\
\mbox{subject to} & {\displaystyle \frac{a_{k}^{2}\|x\|_{2}^{2}}{4\lambda}+a_{k}\langle x,\widehat{\xi}_{i}\rangle + b_{k}\tau  \leq s_{i},  } & \forall i\in[N], k\in\{1,2\}
\end{array}
\right. \label{eqn2:ReformulCVaRCasoP2}
\end{align}
where equation (\ref{eqn1:ReformulCVaRCasoP2}) follows from Theorem \ref{Thm:ReformulacionDROWInterno}. Upon initial examination, the problem obtained may not appear to belong to any known class of optimization problems. However, with some algebraic manipulation, it is possible to express this problem as a linear problem with quadratic constraints. Nonetheless, the matrices describing the quadratic constraints are not semidefinite, which presents a disadvantage. To overcome this, we note that if we fix $\lambda$ and treat it as a constant rather than a variable, the resulting problem associated with that $\lambda$ is once again a linear problem with quadratic constraints. In this case, the matrix defining the quadratic constraints is semidefinite. Therefore, to solve the problem (\ref{eqn2:ReformulCVaRCasoP2}) in this work, we search for the value of $\lambda$ that generates the smallest optimal value in its associated linear problem with quadratic constraints. We achieve this through the use of derivative-free optimization techniques, specifically a method similar to the Nelder-Mead method. Although this may not be the most efficient way to solve this problem, it is observed that the problem obtained by the standard approach (\ref{eqn2:ReformulCVaRCasoP2}) is not trivial to solve.

Now, we introduce the alternative distributional robust approach for this scenario. To this end, it is essential to evaluate the 2-Lipschitz norm $\gamma_{x,F,2}$ of the function $F(x, \xi)$ for this context, which can be readily demonstrated to be $\gamma_{x,F,2}=\|x\|_{2}\max\{|a_{1}|,|a_{2}|\}$. Consequently, the alternative distributional robust approach stated in equation (\ref{CVaRWassertsAlternative}), can be formulated and reformulated as follows:
\begin{align}
    \widehat{J}^{\mathrm{A}}_{N,p,q}(\varepsilon)
    &=  \inf_{x\in\mathcal{X},\tau\in\mathbb{R}}\left( \varepsilon \|x\|_{2}\max_{k\in\{1,2\}}|a_{k}|+\frac{1}{N}\sum_{i=1}^{N}\max_{k\in\{1,2\}}\left(a_{k}\langle a,\widehat{\xi}_{i}\rangle +b_{k}\tau\right) \right) \label{eqn1:ReformulCVaRCasoP2Alter}\\
    &= \left\{
\begin{array}{lll}
{\displaystyle \inf_{x\in\mathcal{X},\tau\in\mathbb{R},\lambda\geq 0,s\in\mathbb{R}^{N}}} & {\displaystyle \lambda \varepsilon +\frac{1}{N}\sum_{i=1}^{N}s_{i}} &\\
\mbox{subject to} & {\displaystyle a_{k}\langle x,\widehat{\xi}_{i}\rangle + b_{k}\tau  \leq s_{i},  } & \forall i\in[N], k\in\{1,2\}\\
& {\displaystyle \|x\|_{2}\max_{k\in\{1,2\}}|a_{k}|\leq \lambda} &
\end{array}
\right. \label{eqn2:ReformulCVaRCasoP2Alter}
\end{align}
where equation (\ref{eqn1:ReformulCVaRCasoP2Alter}) follows from Theorem \ref{Prop:Reformul1MarkovizDROWReformLargaLargaSupportAcot}(b). The problem resulting (\ref{eqn2:ReformulCVaRCasoP2Alter})  can be formulated as a second-order cone program. In comparison with the standard approach presented in (\ref{eqn2:ReformulCVaRCasoP2}), the alternative method seems to generate a computationally easier problem for this specific case.

\begin{figure}[t]  
\centering
 \begin{tabular}{ccc}
    \includegraphics[scale=0.27]{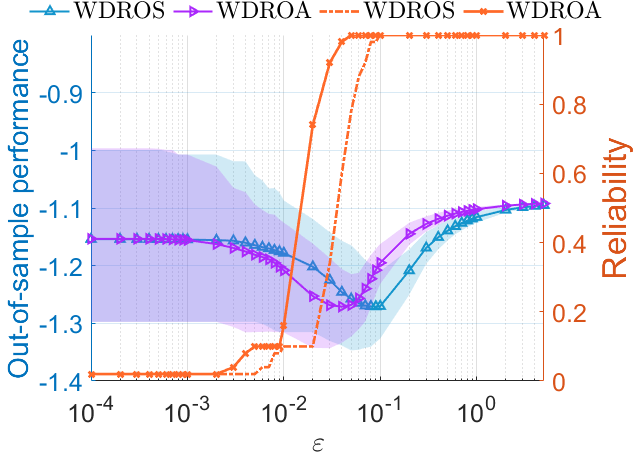} & \includegraphics[scale=0.27]{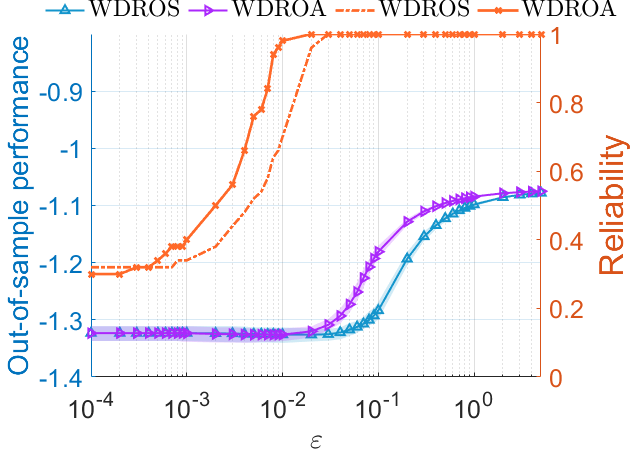}  & \includegraphics[scale=0.27]{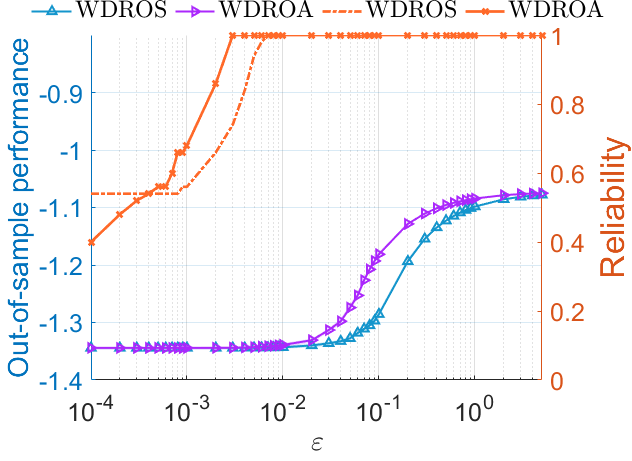}\\
    (a) $N=30$ & (b) $N=300$ & (c)$N=3000$\\
    \includegraphics[scale=0.27]{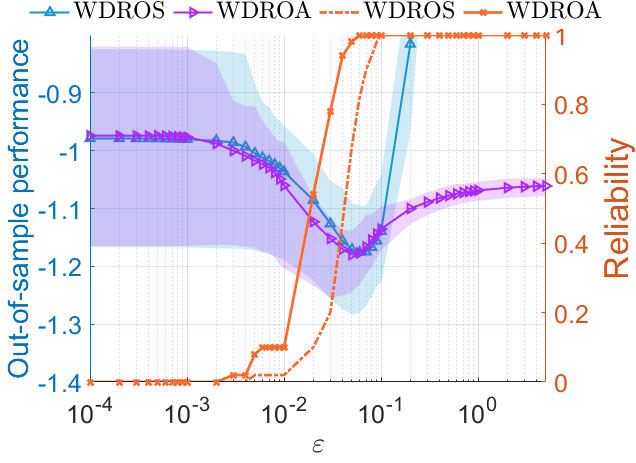}& \includegraphics[scale=0.27]{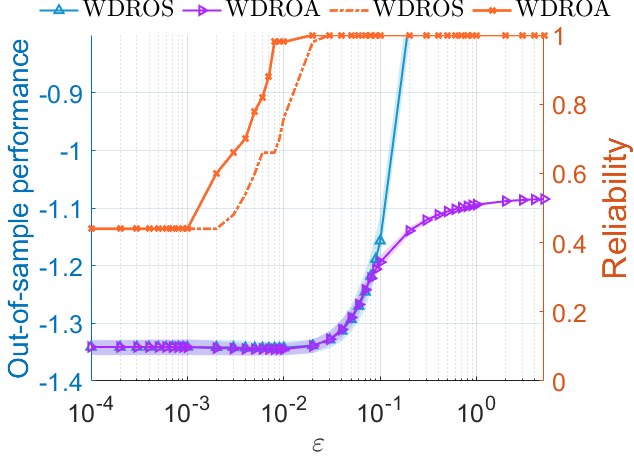}  & \includegraphics[scale=0.27]{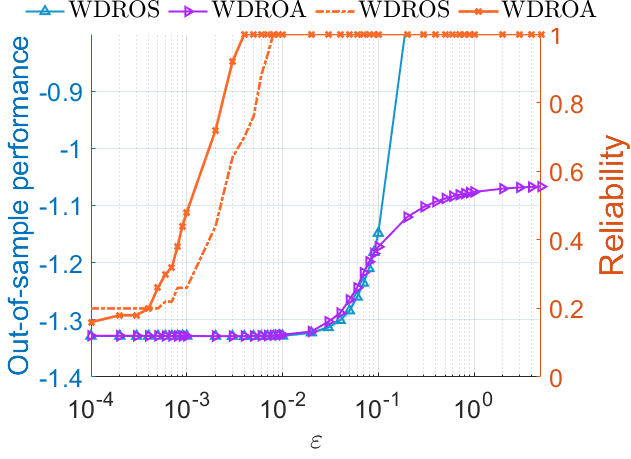}  \\
    (d)$N=30$ & (e)$N=300$  & (f) $N=3000$ \\
    \includegraphics[scale=0.28]{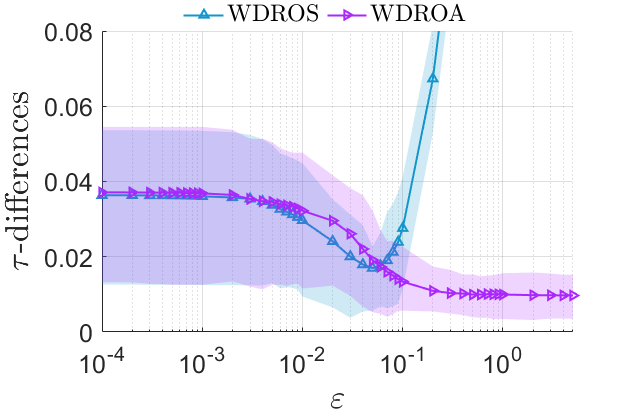} &  \includegraphics[scale=0.28]{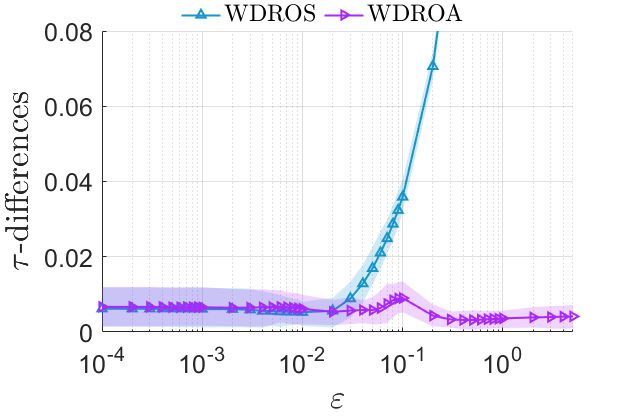}  &  \includegraphics[scale=0.28]{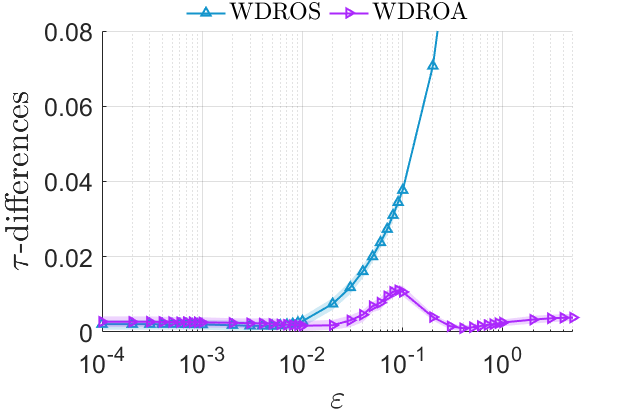} \\
    (g) $N=30$  & (h) $N=300$  & (i) $N=3000$ 
 \end{tabular}
\caption{ (a)-(c) Portfolio-focused out-of-sample performance  (left axis, solid lines, and shaded areas) and
reliability (right axis, solid lines) (d)-(f) Comprehensive out-of-sample performance (left axis, solid lines, and shaded areas) and
reliability (right axis, solid lines). (g)-(i) $\tau$-differences, this refers to the discrepancy between the values of $\widehat{\tau}_{N}^{\mathrm{S}}(\varepsilon)$ and $\widehat{\tau}_{N}^{\mathrm{A}}(\varepsilon)$, and the corresponding value of $\tau$ derived from solving equation (\ref{eqn:DefOutOfSamplePerformSinTau}) for the respective $\widehat{x}_{N}^{\mathrm{S}}(\varepsilon)$ and $\widehat{x}_{N}^{\mathrm{A}}(\varepsilon)$. All are described as a function of the Wasserstein  radius $\varepsilon$ and estimated on the basis of 50 simulations. } \label{fig:MeanCVaRSimulacionesNonLimited}
\end{figure}

Figure \ref{fig:MeanCVaRSimulacionesNonLimited} depicts the findings of the conducted numerical experiments. The results are shown for every $N\in \{30,300,3000\}$, where Figures \ref{fig:MeanCVaRSimulacionesNonLimited}(a)-(c) illustrate the portfolio-focused out-of-sample performance and  reliability. Additionally, Figures \ref{fig:MeanCVaRSimulacionesNonLimited}(d)-(f) showcase comprehensive out-of-sample performance and reliability, while Figures \ref{fig:MeanCVaRSimulacionesNonLimited}(g)-(i) show the dissimilarities between the values of $\tau$ generated by solving (\ref{eqn:DefOutOfSamplePerformSinTau}) for portfolios produced by each approach, measured in absolute value and the values $\widehat{\tau}_{N}^{\mathrm{S}}(\varepsilon)$ and $\widehat{\tau}_{N}^{\mathrm{A}}(\varepsilon)$ respectively. All shaded regions in these figures represent the results of the simulations for each case, these are the tubes between the 20\% and 80\% quantiles, and the solid lines are the averages.

Based on Figures \ref{fig:MeanCVaRSimulacionesNonLimited}(a)-(c), the portfolio-focused out-of-sample performance of the WDROA approach appears to be a contraction of the performance of the WDROS approach as a function of $\varepsilon$. It is not possible to conclude from these experiments which approach is superior in terms of this performance. However, both approaches exhibit comparable portfolio-focused out-of-sample performance. For instance, in the context of Figure \ref{fig:MeanCVaRSimulacionesNonLimited}(a), we observe that the minimum average portfolio-focused out-of-sample performance achieved with the standard approach is nearly identical to the minimum average for this performance achieved with the alternative approach. This observation is significant because, in practice, one aims to estimate the $\varepsilon$ value that corresponds to the minimum portfolio-focused out-of-sample performance, which can be achieved using methods such as the Holdout or $k$-fold cross-validation described in \cite{Kuhn2018}. Regardless of the approach used, employing the same method will result in $\varepsilon$  values that yield portfolios with similar portfolio-focused out-of-sample performance, although the $\varepsilon$ values will not be identical. Thus, the robust approach guarantees portfolio-focused out-of-sample performance that is similar to that corresponding to the standard approach but with an optimization problem (\ref{eqn2:ReformulCVaRCasoP2Alter}) that appears to be less demanding than that of the standard approach (\ref{eqn2:ReformulCVaRCasoP2}). 

Regarding comprehensive out-of-sample performance, the results are displayed in Figure \ref{fig:MeanCVaRSimulacionesNonLimited}(d)-(f). The figures illustrate that as $\epsilon$ increases, the WDROS method yields significantly high values, in contrast to WDROA. The findings concerning comprehensive out-of-sample performance have led to an important observation. Let $\tau^{*}(\widehat{x}_{N}^{\mathrm{A}}(\varepsilon))$ denote the optimal solution of (\ref{eqn:DefOutOfSamplePerformSinTau}) when $\widehat{x}_{N}=\widehat{x}_{N}^{\mathrm{A}}(\varepsilon)$. Note that, by Theorem 1 of \cite{Rockafellar2000}, if the set of optimal solutions of (\ref{eqn:DefOutOfSamplePerformSinTau}) is a single point, then $\tau^{*}(\widehat{x}_{N}^{\mathrm{A}}(\varepsilon))$ is the VaR of portfolio $\widehat{x}_{N}^{\mathrm{A}}(\varepsilon)$. In other words, $\tau^{*}(\widehat{x}_{N}^{\mathrm{A}}(\varepsilon))=\mathrm{VaR}_{\alpha,\xi\sim\mathbb{P}}\left(-\langle \widehat{x}_{N}^{\mathrm{A}}(\varepsilon),\xi\rangle\right)$. The same holds for  $\widehat{x}_{N}^{\mathrm{S}}(\varepsilon)$. Therefore, it is to be expected that both $\widehat{\tau}_{N}^{\mathrm{S}}(\varepsilon)$ and $\widehat{\tau}_{N}^{\mathrm{A}}(\varepsilon)$ provide information about $\tau^{*}(\widehat{x}_{N}^{\mathrm{S}}(\varepsilon))$ and $\tau^{*}(\widehat{x}_{N}^{\mathrm{A}}(\varepsilon))$ respectively. Since $\mathbb{P}$ is generally unknown in practice, $\tau^{*}(\widehat{x}_{N}^{\mathrm{A}}(\varepsilon))$ is also unknown. However, as $\tau^{*}(\widehat{x}_{N}^{\mathrm{A}}(\varepsilon))$ is employed to determine portfolio-focused out-of-sample performance, it follows from Figures \ref{fig:MeanCVaRSimulacionesNonLimited}(a)-(c) and \ref{fig:MeanCVaRSimulacionesNonLimited}(d)-(f) that the portfolio-focused and comprehensive out-of-sample performances tend to exhibit similar shapes and scales for the alternative approach, implying that $\widehat{\tau}_{N}^{\mathrm{A}}(\varepsilon)$ approximates $\tau^{*}(\widehat{x}_{N}^{\mathrm{A}}(\varepsilon))$. The latter is illustrated in Figures \ref{fig:MeanCVaRSimulacionesNonLimited}(g)-(i). This, however, does not hold for the standard approach because significant differences between $\tau^{*}(\widehat{x}_{N}^{\mathrm{S}}(\varepsilon))$ and $\widehat{\tau}_{N}^{\mathrm{S}}(\varepsilon)$ emerge when $\varepsilon$ is large. Therefore, if we require an approach that generates acceptable estimators of the VaR of the portfolio, the alternative approach tends to be a better choice.

In practical applications, it is desirable that the optimal solutions obtained from any of the approaches be close to the optimal solution of the problem (\ref{eqn:Mean-CVaROriginalReformu}). To achieve this goal, an appropriate value of $\varepsilon$ needs to be determined. The optimal $\varepsilon$ is the one that minimizes the average comprehensive out-of-sample performance (solid lines in Figures \ref{fig:MeanCVaRSimulacionesNonLimited}(d)-(f)). For instance, in Figure \ref{fig:MeanCVaRSimulacionesNonLimited}(d), the optimal $\varepsilon$ value for the two compared approaches is around $6\times10^{-2}$. However, that $\varepsilon$ is unknown in practice, and thus, it needs to be estimated from the sample data. Commonly used methods for this purpose include Holdout and $k$-fold cross-validation. However, these methods have a margin of error and when the sample size is small, the margin of error of the estimation can be high. If the solutions of both approaches fall on the right end of the interval generated by the margin of error, the standard approach solutions $\widehat{x}_{N}^{\mathrm{S}}(\varepsilon)$ and $\widehat{\tau}_{N}^{\mathrm{S}}(\varepsilon)$ would have a higher comprehensive out-of-sample performance than the alternative approach solutions for that $\varepsilon$. Based on this type of performance, it can be concluded that the alternative approach outperforms the standard approach, particularly for larger values of $\epsilon$.

In terms of reliability, the desired outcome is that the reliability approaches 1 as $\varepsilon$ increases, for both portfolio-focused and comprehensive reliability. Figures \ref{fig:MeanCVaRSimulacionesNonLimited}(a)-(c) and (d)-(i) confirm that this expectation is met. Moreover, it is observed that, again, the alternative approach can be seen as a contraction of the standard one.

%the standard approach attains higher reliability values initially. Achieving high reliability for any approach is critical because it guarantees that the optimal value of the optimization problem connected to the approach exceeds the optimal value of the original problem   with a high probability. The original problem corresponds to equation (\ref{eqn:Mean-CVaROriginal}) for the portfolio-focused case and equation (\ref{eqn:Mean-CVaROriginalReformu}) for the comprehensive case.

\begin{figure}[t]  
\centering
 \begin{tabular}{ccc}
 %   \includegraphics[scale=0.28]{PortfoliosOurN30.png} &  \includegraphics[scale=0.28]{PortfoliosOurN300.png}  &  \includegraphics[scale=0.28]{PortfoliosOurN3000.png} \\
  %  (a) $N=30$ WDROA & (b) $N=300$ WDROA & (c) $N=3000$ WDROA\\
 %    \includegraphics[scale=0.28]{PortfoliosTradN300.png}& \includegraphics[scale=0.28]{PortfoliosTradN300.png} & \includegraphics[scale=0.28]{PortfoliosTradN3000.png} \\
  %   (d) $N=30$ WDROS & (e) $N=300$ WDROS & (f) $N=3000$ WDROS \\
    \includegraphics[scale=0.27]{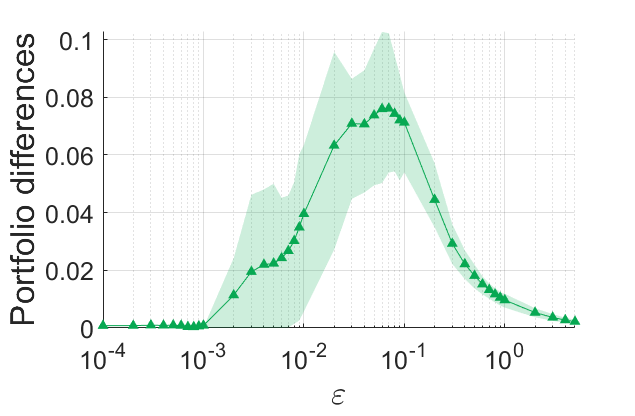} & \includegraphics[scale=0.27]{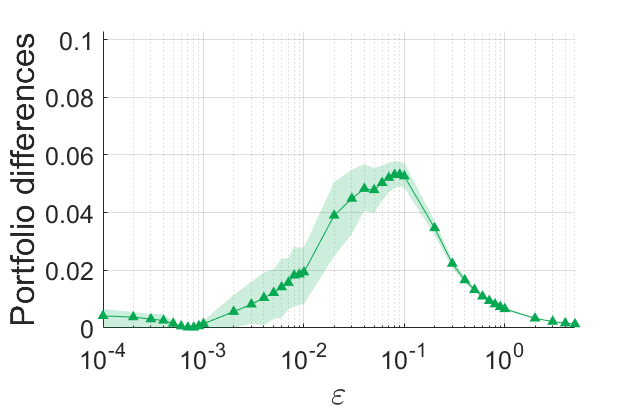}  & \includegraphics[scale=0.27]{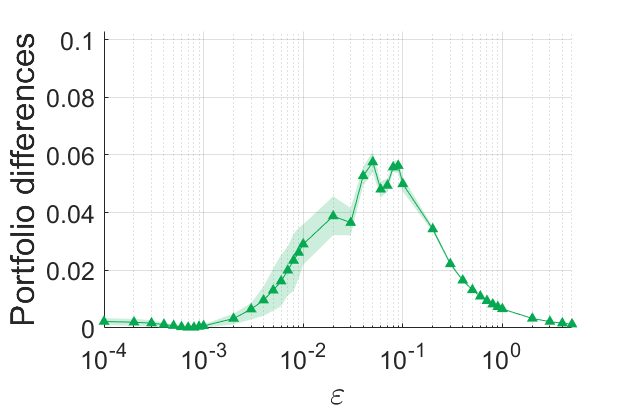} \\
     (a)$N=30$ & (b) $N=300$ & (c) $N=3000$
     \end{tabular}
\caption{ 
%(a)-(f) Optimal composition of the portfolio as a function of the Wasserstein radius $\varepsilon$ for each approach (WDROS and WDROA), obtained by averaging over 50 simulations. The weights assigned to each asset are depicted in ascending order, where the weight assigned to asset 1 is displayed at the bottom (dark blue area) and that assigned to asset 10 at the top (dark red area). 
(a)-(c) Portfolio differences in $\infty$-norm obtained by averaging over 50 simulations.} \label{fig:MeanCVaRSimulacionesPortfolios}
\end{figure}

%Figure \ref{fig:MeanCVaRSimulacionesPortfolios} depicts the behavior of portfolios for varying values of $\varepsilon$ and $N$. The figure displays the average of portfolios obtained from 50 simulations, differentiated by a color scale. Specifically, Figures \ref{fig:MeanCVaRSimulacionesPortfolios}(a)-(c) represent the standard approach, while Figures \ref{fig:MeanCVaRSimulacionesPortfolios}(d)-(f) represent the alternative approach proposed in this study. For both approaches, the average portfolio converges to the equally weighted portfolio $x^{*}=e/m$, which aligns with Proposition 7.2 in \cite{Kuhn2018}. Furthermore, 
Figures \ref{fig:MeanCVaRSimulacionesPortfolios}(a)-(c) show the contrast between the portfolios generated by the two approaches studied. Given a sample, the plot shows the $\|\widehat{x}_{N}^{\mathrm{S}}(\varepsilon)-\widehat{x}_{N}^{\mathrm{A}}(\varepsilon) \|_{\infty}$ value for each $\varepsilon$, repeated 50 times, resulting in the green shaded region, which is the tube between the 20\% and 80\% quantiles. The solid line denotes the average. Notably, a region of $\varepsilon$ values exists where the difference between the portfolios increases. As $N$ becomes larger, the difference stabilizes but does not become negligible. For both approaches, the average portfolio converges to the equally weighted portfolio $x^{*}=e/m$, which aligns with Proposition 7.2 in \cite{Kuhn2018}.
%Although simulations have produced promising results, there is still a need to explore the theoretical underpinnings of these findings. Therefore, further research is necessary to demonstrate the validity of these results from a theoretical perspective.

\subsection{Limited-loss assets}

In this part, the objective is to observe that there are cases where although theoretically it is known that the two approaches are equivalent, sometimes the formulation does not seem to be equivalent. Here, we adopt the $1$-Wasserstein distance and cost function $\mathbf{d}=\|\cdot\|_{1}$. Additionally, we consider $\Xi=\left\{\xi\in\mathbb{R}^{m}\: : \:\xi\geq -1\right\}$. Corollary 5.1 of \cite{Kuhn2018}, which follows from Theorem \ref{Thm:ReformulacionDROWInterno}, reveals that the standard distributionally robust formulation WDROS that employs this Wasserstein distance is equivalent to the following optimization problem:
\begin{equation} \label{eqn2:ReformulCVaRCasoP1}
\widehat{J}^{\mathrm{S}}_{N,p,q}(\varepsilon) = \left\{
\begin{array}{lll}
{\displaystyle \inf_{x\in\mathcal{X},\tau\in\mathbb{R},\lambda\geq 0,s\in\mathbb{R}^{N},\gamma}} & {\displaystyle \lambda \varepsilon +\frac{1}{N}\sum_{i=1}^{N}s_{i}} &\\
\mbox{subject to} & {\displaystyle \langle \gamma_{i,k},e+\widehat{\xi}_{i} \rangle+a_{k}\langle x,\widehat{\xi}_{i}\rangle + b_{k}\tau  \leq s_{i},  } & \forall i\in[N], k\in\{1,2\} \\
& \|\gamma_{i,k}+a_{k} x  \|_{\infty} \leq \lambda, & \forall i\in[N], k\in\{1,2\} \\
 & \gamma_{i,k}\geq 0 . & \forall i\in[N], k\in\{1,2\} 
\end{array}
\right. 
\end{equation}
The problem resulting (\ref{eqn2:ReformulCVaRCasoP1})  is a linear program with $N(1+2m)+m+2$ variables.
\vspace{-0.2cm}
\begin{equation*} %Este parrafo lo introduje de esta forma porque las ecuaciones tienen problemas.
  \parbox{4.8in}{%
    \qquad Now, we introduce the alternative distributional robust approach for this scenario. To this end, it is essential to evaluate the $\infty$-Lipschitz norm $\gamma_{x,F,\infty}$ of the function $F(x, \xi)$ for this context, which can be readily demonstrated to be $\gamma_{x,F,\infty}=\|x\|_{\infty} \max\limits_{k=1,2}|a_{k}|$. This last supremum is reached at $\infty$, so by Corollary 2 of \cite{GaoChenKleywegt2017}, the two approaches are equivalent. Consequently, since it is not hard to show that $\inf\limits_{x\in\mathcal{X},\tau\in\mathbb{R}}\sup\limits_{\xi\in \Xi}\left(\max\limits_{k=1,2}(a_{k}\langle x\xi\rangle +b_{k}\tau)\right)=\infty$, the alternative distributional robust approach WDROA can be formulated and reformulated as follows:%
}
\end{equation*}

\begin{align}
    \widehat{J}^{\mathrm{A}}_{N,p,q}(\varepsilon) &= \inf_{x\in\mathcal{X},\tau\in\mathbb{R}}\left(\varepsilon \|x\|_{\infty}\max_{k\in\{1,2\}}|a_{k}|+\frac{1}{N}\sum_{i=1}^{N}\max_{k\in\{1,2\}}\left(a_{k}\langle a,\widehat{\xi}_{i}\rangle +b_{k}\tau\right) \right) \nonumber \\
    &= \left\{
\begin{array}{lll}
{\displaystyle \inf_{x\in\mathcal{X},\tau\in\mathbb{R},\lambda\geq 0,s\in\mathbb{R}^{N}}} & {\displaystyle \lambda \varepsilon +\frac{1}{N}\sum_{i=1}^{N}s_{i}} &\\
\mbox{subject to} & {\displaystyle a_{k}\langle x,\widehat{\xi}_{i}\rangle + b_{k}\tau  \leq s_{i},  } & \forall i\in[N], k\in\{1,2\},\\
& {\displaystyle \|x\|_{\infty}\max_{k\in\{1,2\}}|a_{k}|\leq \lambda} &
\end{array}
\right. \label{eqn2:ReformulCVaRCasoP1Alter}
\end{align}
 The problem resulting (\ref{eqn2:ReformulCVaRCasoP1Alter})  can be formulated as a linear program with $m+N+2$. In comparison with the standard approach presented in (\ref{eqn2:ReformulCVaRCasoP1}), the alternative method seems to have fewer variables for this specific case. Furthermore, if the optimal solutions of (\ref{eqn2:ReformulCVaRCasoP1}) linked to  $\gamma_{i,k}$ variables were zero, both problems would be equal. However, it was observed by means of numerical experiments that not all of the optimal solutions associated with  $\gamma_{i,k}$ variables were zero. Specifically, between $\gamma_{i,1}$ and $\gamma_{i,2}$ always one is equal to zero and the other is not. Thus, it remains unclear how to transform one formulation into the other. This suggests that the alternative approach results in a formulation that computationally may be more convenient than the one derived from the standard approach.

\textbf{Acknowledgments}\\
This work was supported by Fondo de Investigaciones de la Facultad de Ciencias de la Universidad de los Andes INV-2021-128-2307 and INV-2021-126-2273.

\bibliographystyle{splncs04}
\bibliography{refs}
%

%**********************************************

\appendix

\section{Proofs of Lemmas and Theorems}

We present proof of the results presented in this work. Section \ref{Apendice:PruebasLema} presents the proof of Lemma \ref{Lemma:JustificacionBola}. Section \ref{Apendice:RestriccionesVE} explores worst-case expectation problems with expected value constraints and the corresponding dual formulation. Theorem \ref{Thm:DualidadFuerteDROWconRestric} is an important result on its own. Section \ref{Apendice:DRVarianza} presents a distributionally robust estimation of the variance under known mean, a key result to prove Theorem \ref{Prop:Reformul1MarkovizDROWReformLargaLarga}. Section \ref{Sec:Appendix:ProofThmMeanVarianze} shows the proofs of Theorems \ref{Prop:Reformul1MarkovizDROWReformLargaLargaSupportAcot} and \ref{Prop:Reformul1MarkovizDROWReformLargaLarga} and its corollaries. 

\subsection{Proof of Lemma \ref{Lemma:JustificacionBola}.}\label{Apendice:PruebasLema}

\proof[Lemma \ref{Lemma:JustificacionBola}] Let  $\widetilde{\xi}_{1},\ldots,\widetilde{\xi}_{M}$ be another sample of $\xi$, then let   $\widetilde{\mathbb{P}}_{M}$ be the empirical distribution generated by this sample. This last sample of $\xi$ induces the sample $\widetilde{\zeta}^{x,F}_{1},\ldots,\widetilde{\zeta}^{x,F}_{M}$ of $\zeta^{x,F}$ given by  $\widetilde{\zeta}^{x,F}_{i}:=F\left( x, \widetilde{\xi}_{i}\right)$, so we consider $\widetilde{\mathbb{P}}^{x,F}_{M}$, the empirical distribution generated by the sample  $\left\{\widetilde{\zeta}^{x,F}_{i}\right\}_{i=1}^{M}$. Because $\widetilde{\mathbb{P}}_{M} \rightarrow \mathbb{P}$ and $ \widetilde{\mathbb{P}}^{x,F}_{M} \rightarrow \mathbb{P}^{x,F}$ weakly as $M$ goes to $\infty$, by Corollary 6.11 en \cite{Villani2008} we have that
\begin{equation} \label{Eq:ConvergenciaEnWassesrtenEmiric}
W_{p}\left( \widehat{\mathbb{P}}_{N} ,\widetilde{\mathbb{P}}_{M}\right)\overset{{\scriptstyle M\rightarrow \infty}}{\longrightarrow}W_{p}\left( \widehat{\mathbb{P}}_{N} ,\mathbb{P}\right) \quad \mbox{ and }\quad W_{p}\left( \widehat{\mathbb{P}}^{x}_{N} ,\widetilde{\mathbb{P}}^{x,F}_{M}\right)\overset{{\scriptstyle M\rightarrow \infty}}{\longrightarrow}W_{p}\left( \widehat{\mathbb{P}}_{N}^{x,F} ,\mathbb{P}^{x,F}\right). 
\end{equation}
Additionally, for each $M$ we get that
\begin{align*}
W_{p}^{p}\left( \widehat{\mathbb{P}}^{x,F}_{N} ,\widetilde{\mathbb{P}}^{x,F}_{M}\right) &= \inf\left\{\sum_{i=1}^{N}\sum_{j=1}^{M} \lambda_{i,j}\left|\widehat{\zeta}^{x,F}_{i}-\widetilde{\zeta}^{x,F}_{j}  \right|^{p} \:\left|\: \begin{array}{l} \sum_{i=1}^{N}\lambda_{i,j}=\frac{1}{M},\\ \sum_{j=1}^{M}\lambda_{i,j}=\frac{1}{N},\\ \lambda_{i,j}\geq 0, \\ i=1,\ldots,N,\\ j=1,\ldots,M \end{array}  \right.\right\} \\
&= \inf\left\{\sum_{i=1}^{N}\sum_{j=1}^{M} \lambda_{i,j}\left|F\left( x, \widehat{\xi}_{i}\right)-F\left( x, \widetilde{\xi}_{j}\right)   \right|^{p} \:\left|\: \begin{array}{l} \sum_{i=1}^{N}\lambda_{i,j}=\frac{1}{M},\\ \sum_{j=1}^{M}\lambda_{i,j}=\frac{1}{N},\\ \lambda_{i,j}\geq 0, \\ i=1,\ldots,N,\\ j=1,\ldots,M \end{array}  \right.\right\}  \\
&\leq  \inf\left\{\sum_{i=1}^{N}\sum_{j=1}^{M} \lambda_{i,j} \gamma_{x,F}^{p}\left\| \widehat{\xi}_{i}- \widetilde{\xi}_{j}  \right\|^{p} \:\left| \: \begin{array}{l} \sum_{i=1}^{N}\lambda_{i,j}=\frac{1}{M},\\ \sum_{j=1}^{M}\lambda_{i,j}=\frac{1}{N},\\ \lambda_{i,j}\geq 0, \\ i=1,\ldots,N,\\ j=1,\ldots,M \end{array}  \right.\right\} \:\: \begin{array}{l}\mbox{by Assumption \ref{AssumptionPrincipal}}\\ \mbox{}\end{array}\\
&= \gamma_{x,F}^{p}W_{p}^{p}\left( \widehat{\mathbb{P}}_{N} ,\widetilde{\mathbb{P}}_{M}\right).
\end{align*} 
Therefore, by (\ref{Eq:ConvergenciaEnWassesrtenEmiric}) we conclude 
$$W_{p}^{p}\left( \widehat{\mathbb{P}}_{N}^{x,F} ,\mathbb{P}^{x,F}\right)\leq \gamma_{x,F}^{p} W_{p}^{p}\left( \widehat{\mathbb{P}}_{N} ,\mathbb{P}\right).$$
\qed
\endproof

\subsection{Duality of worst-case expectation problems with expected value constraints}\label{Apendice:RestriccionesVE}

When we refer to worst-case expectation  problems with restrictions on the expected value, we are referring to problems of the form

\begin{equation}\label{Eqn:VraicionDROW}
\left\{
\begin{array}{lll}
{\displaystyle \sup_{\mathbb{Q}\in \mathcal{B}_{\varepsilon}(\widehat{\mathbb{P}}_{N})} } & \mathbb{E}_{\mathbb{Q}}\left[ h(\xi) \right] \\
\mbox{subject to  } & \mathbb{E}_{\mathbb{Q}}\left[g_{i}(\xi)\right]=b_{i}. & \forall \: i=1,\ldots,k.
\end{array}
\right.
\end{equation}
where $b_{i}\in\mathbb{R}$ and $g_{1},\ldots,g_{k}$ are integrable functions respect to each measure in $\mathcal{B}_{\varepsilon}(\widehat{\mathbb{P}}_{N})$. This problem is important for the following section and to prove Theorem \ref{Thm:EsimaVarConMediaConoc}.

\begin{theorem} \label{Thm:DualidadFuerteDROWconRestric}
Assume that the optimal value of the problem (\ref{Eqn:VraicionDROW}) is finite and that any of the following conditions are satisfied
\begin{enumerate}
\item[i)] The point $(b_{1},\ldots,b_{k},1)$ is  a interior point of the set
$$\left\{\left.\lambda\left( \int g_{1}(\xi)\mathbb{Q}(d\xi),\ldots, \int g_{k}(\xi)\mathbb{Q}(d\xi),1\right) \:\right|\: \lambda>0,\:\mathbb{Q}\in \mathcal{B}_{\varepsilon}(\widehat{\mathbb{P}}_{N}) \right\},$$
\item[ii)] The set of optimal distributions of  (\ref{Eqn:VraicionDROW}) is not empty and bounded.
\end{enumerate}
Then (\ref{Eqn:VraicionDROW}) satisfies strong duality, that is, the optimal value of (\ref{Eqn:VraicionDROW}) is equal to
$$\inf_{a_{1},\ldots, a_{k}} \left\{ \sum_{i=1}^{k}a_{i}b_{i}+\sup_{\mathbb{Q}\in \mathcal{B}_{\varepsilon}(\widehat{\mathbb{P}}_{N})}\int_{\Xi}\left(h(\xi)-\sum_{i=1}^{k}a_{i}g_{i}(\xi)\right)\mathbb{Q}(d\xi) \right\}.$$
\end{theorem}
The proof of this theorem consists of reformulating (\ref{Eqn:VraicionDROW})  as a linear conic problem and in that context using the results of \cite{Shapiro2001} to show that strong duality is satisfied.

Using Theorem \ref{Thm:ReformulacionDROWInterno} we can further reformulate (\ref{Eqn:VraicionDROW}) as a semi-infinite optimization problem.

\begin{corollary} \label{Corol:DualidadFuerteDROWconRestric}
Suppose that the function $F_{a}(\xi):=h(\xi)-\sum_{i=1}^{k}a_{i}(g_{i}(\xi)-b_{i})$ satisfies the hypotheses of the Theorem \ref{Thm:ReformulacionDROWInterno} for all $a\in\mathbb{R}^{k}$, and  satisfies any of the conditions i) and ii) of Theorem \ref{Thm:DualidadFuerteDROWconRestric}, then   the problem (\ref{Eqn:VraicionDROW}) can be rewritten as

\begin{equation}
\left\{\begin{array}{lll}
{\displaystyle \inf_{a_{1},\ldots,a_{k},\lambda} }& {\displaystyle \lambda\varepsilon^{p}+  \frac{1}{N}\sum_{i=1}^{N}s_{i}  }\\
\mbox{subject to  } & {\displaystyle  \sup_{\xi\in\Xi}\left(h(\xi)-\sum_{i=1}^{k}a_{i}(g_{i}(\xi)-b_{i} )-\lambda d^{p}(\xi,\widehat{\xi}_{i}) \right)\leq s_{i}  } & \forall\: i=1,\ldots,N,\\
& \lambda\geq 0.
\end{array}\right.
\end{equation}
\end{corollary}

\subsection{Distributionally robust estimation of the variance of a random variable with known mean}\label{Apendice:DRVarianza}

In this part, we will formulate a robust distributional version of the problem of estimating the variance of a random variable with a known mean, and we will demonstrate that the obtained optimization problem admits an explicit solution.

Let $\zeta$ be a random variable with unknown distribution $\mathbb{P}$ with support $\Xi\subseteq\mathbb{R}$, we assume that the expected value of $\zeta$ is known, specifically, we assume that $\mathbb{E}_{\mathbb{P}}[\zeta]=\eta$. Also, we consider a sample $\widehat{\zeta}_{1},\ldots,\widehat{\zeta}_{N}$ of $\zeta$. Let $\widehat{\mathbb{P}}_{N}:=\frac{1}{N}\sum_{i=1}^{N}\delta_{\widehat{\zeta}_{i}}$ be the empirical distribution, and denote by $\bar{\zeta}:=\mathbb{E}_{\widehat{\mathbb{P}}_{N}}[\zeta]=\frac{1}{N}\sum_{i=1}^{N}\widehat{\zeta}_{i}$ and $\widehat{\sigma}^2:=\mathbb{E}_{\widehat{\mathbb{P}}_{N}}[(\zeta-\eta)^2]=\frac{1}{N}\sum_{i=1}^{N}(\widehat{\zeta}_{i}-\eta)^2$, the empirical mean and variance respectively. Let $\mathcal{B}_{\varepsilon}(\widehat{\mathbb{P}}_{N})$ be the $2$-Wasserstein ball with center in $\widehat{\mathbb{P}}_{N}$ and radius $\varepsilon$. We call the following problem a distributionally robust estimate of the variance of $\zeta$:

\begin{align}\left\{\begin{array}{ll}
{\displaystyle \sup_{\mathbb{Q}\in\mathcal{B}_{\varepsilon}(\widehat{\mathbb{P}}_{N}) } } & \mathbb{E}_{\mathbb{Q}}\left[\left(\zeta -\eta \right)^{2}  \right]\\
\mbox{subject to  } & \mathbb{E}_{\mathbb{Q}}\left[\zeta\right]= \eta.
\end{array}  \right.\label{VarianConMeanDROW}
\end{align}  
However, for some values of $\varepsilon$, this problem may be not feasible as the following result shows. 

\begin{proposition}\label{Prop:EpsilonFactible}
If $\varepsilon < \left|\eta -\bar{\zeta} \right|$ then
$$\mathcal{B}_{\varepsilon}(\widehat{\mathbb{P}}_{N}) \cap \left\{ \mathbb{Q}\in\mathcal{P}(\mathbb{R})\:\left|\:\mathbb{E}_{\mathbb{Q} }[\zeta]=\eta\right.\right\}=\emptyset.$$
\end{proposition}

\proof
Let $\mathbb{Q}\in \mathcal{B}_{\varepsilon}(\widehat{\mathbb{P}}_{N})$, we must show that $\mathbb{E}_{\mathbb{Q}}[\zeta]\neq \eta$. Indeed, by Observation 6.6 in \cite{Villani2008} we know that if $p\leq q$ then $W_{p}\leq W_{q}$. In particular, we have that $W_{1}\leq W_{2}$ and this implies that 
\begin{equation}\label{Eqn:DeswualdadW1w2}
W_{1}\left(\mathbb{Q},\widehat{\mathbb{P}}_{N} \right)\leq W_{2}\left(\mathbb{Q},\widehat{\mathbb{P}}_{N}\right)\leq \varepsilon < \left|\eta -\bar{\zeta} \right|.
\end{equation}
Therefore, defining $\mathcal{S}(\mathbb{Q},\widehat{\mathbb{P}}_{N})$ as the set of couplings between   $\mathbb{Q}$ and $\widehat{\mathbb{P}}_{N}$, there exists $\Pi\in\mathcal{S}(\mathbb{Q},\widehat{\mathbb{P}}_{N})$ such that
$$\int_{\Xi\times\Xi}|\zeta-\delta|\Pi(d\xi,d\zeta)<\left|\eta -\bar{\zeta} \right|.$$
We also have
$$
\int_{\Xi\times\Xi}\zeta\Pi(d\zeta,d\delta)=\int_{\Xi}\zeta\mathbb{Q}(d\zeta)=\mathbb{E}_{\mathbb{Q}}[\zeta] \quad\mbox{ and } \quad  \int_{\Xi\times\Xi}\delta\Pi(d\zeta,d\delta)=\int_{\Xi}\delta\widehat{\mathbb{P}}_{N}(d\delta)=\mathbb{E}_{\widehat{\mathbb{P}}_{N}}[\delta].
$$
Then
$$
\left|\mathbb{E}_{\mathbb{Q}}[\zeta]-\bar{\zeta}\right|= \left|\int_{\Xi\times\Xi}(\zeta-\delta)\Pi(d\xi,d\zeta) \right| \leq \int_{\Xi\times\Xi}|\zeta-\delta|\Pi(d\xi,d\zeta).   
$$
In consequence, we obtain
$$\left|\mathbb{E}_{\mathbb{Q}}[\zeta]-\bar{\zeta}\right| < \left|\eta -\bar{\zeta}\right|.
$$
From the above and the inverse triangular inequality follows
$$\left|\eta-\mathbb{E}_{\mathbb{Q}}[\zeta] \right|=\left|\eta-\bar{\zeta}- \left(\mathbb{E}_{\mathbb{Q}}[\zeta]-\bar{\zeta} \right) \right|\geq  \left|\left|\eta -\bar{\zeta} \right|-\left|\mathbb{E}_{\mathbb{Q}}[\zeta]-\bar{\zeta} \right|\right|>0.$$
This allows us to conclude that $\mathbb{E}_{\mathbb{Q}}[\zeta]\neq \eta $. 
\endproof

The following theorem establishes an explicit expression for the optimal value of the optimization problem to the right of (\ref{VarianConMeanDROW}).

\begin{theorem}\label{Thm:EsimaVarConMediaConoc}
Let $\varepsilon>0$  with  $\varepsilon^{2}\geq(\bar{\zeta}-\eta)^2$, and $\Xi=\mathbb{R}$. Then, the optimal value of \eqref{VarianConMeanDROW} is equal to  
$$\left(\sqrt{\widehat{\sigma}^2 -\left(\bar{\zeta}-\eta\right)^{2}  }+\sqrt{\varepsilon^{2}- \left(\bar{\zeta}-\eta\right)^{2} } \right)^{2}.$$
\end{theorem}
\proof 
By Theorem \ref{Thm:DualidadFuerteDROWconRestric} we have that (\ref{VarianConMeanDROW}) satisfies strong duality and its optimal value is equal to
\begin{equation}\label{Eq:Dualidad1}
\inf_{\beta}\sup_{ \mathbb{Q}\in\mathcal{B}_{\varepsilon}(\widehat{\mathbb{P}}_{N}) }\mathbb{E}_{\mathbb{Q}}\left[ \left(\zeta -\eta \right)^{2}-\beta\zeta +\beta\eta  \right].    
\end{equation}
Note that  $\zeta\mapsto\left(\zeta -\eta \right)^{2}-\beta\zeta +\beta\eta $ satisfies the hypotheses of Corollary \ref{Corol:DualidadFuerteDROWconRestric}, therefore this last formulation is equivalent to the semi-infinite optimization program
\begin{align}
\left\{\begin{array}{lll}
{\displaystyle \inf_{\beta,\: \lambda,\: s_{i}} }& {\displaystyle \lambda\varepsilon^{2}+\frac{1}{N}\sum_{i=1}^{N}s_{i}}\\
\mbox{subject to} & {\displaystyle  \sup_{\zeta\in\Xi}\left((\zeta-\eta)^{2}-\beta\zeta+\beta\eta -\lambda \left|\zeta-\widehat{\zeta}_{i} \right|^{2}  \right) \leq s_{i} } & \forall i=1,\ldots,N\\
 & \lambda\geq 0.
\end{array}\right. \label{MarkovizDROWInternoReformSeimiInf}
\end{align} 
If $\lambda<1 $, then $\lambda$ is not a optimal value of (\ref{MarkovizDROWInternoReformSeimiInf}) because, in this case, the set $$\left\{(\zeta-\eta)^{2}-\beta\zeta+\beta\eta -\lambda \left|\zeta-\widehat{\zeta}_{i} \right|^{2}\:|\: \zeta\in\Xi  \right\}$$ is not bounded.  On the other hand, if $\lambda\geq1$, then   $$\sup_{\zeta\in\mathbb{R}}\left((\zeta-\eta)^{2}-\beta\zeta+\beta\eta -\lambda \left|\zeta-\widehat{\zeta}_{i} \right|^{2}  \right)$$    can be calculated explicitly because the function is a concave quadratic polynomial. The unique maximum is attained at $\widehat{\varphi}_{i}=\frac{2\eta+\beta-2\lambda \widehat{\zeta}_{i}}{2(1-\lambda)}$ and (\ref{MarkovizDROWInternoReformSeimiInf}) is equivalent to 
\begin{align}\notag
&\left\{\begin{array}{lll}
{\displaystyle \inf_{\beta,\: \lambda,\: s_{i}} }& {\displaystyle \lambda\varepsilon^{2}+\frac{1}{N}\sum_{i=1}^{N}s_{i}}\\
\mbox{subject to  } & {\displaystyle  \frac{\beta^{2}}{4(\lambda-1)}+\frac{\lambda}{\lambda-1}\left(\beta(\eta-\widehat{\zeta}_{i}) +(\eta-\widehat{\zeta}_{i})^{2}  \right) \leq s_{i} } & \forall i=1,\ldots,N\\
 & \lambda\geq 1.
\end{array}\right.
\\\notag
=&\left\{\begin{array}{ll}
{\displaystyle \inf_{\lambda,\:\beta} } & {\displaystyle \lambda\varepsilon^{2}+ \frac{\beta^{2}}{4(\lambda-1)}+\frac{\lambda}{\lambda-1}\left( \frac{\beta}{N}\sum_{i=1}^{N}(\eta -\widehat{\zeta}_{i}) +\frac{1}{N}\sum_{i=1}^{N}(\eta -\widehat{\zeta}_{i})^{2}   \right) } \\
\mbox{subject to  } & \lambda\geq 1.
\end{array}\right.\\
=&\left\{\begin{array}{ll}
{\displaystyle \inf_{\lambda,\:\beta} } & {\displaystyle \lambda\varepsilon^{2}+ \frac{\beta^{2}}{4(\lambda-1)}+\frac{\lambda}{\lambda-1}\left(\beta(\eta-\bar{\zeta}) +\widehat{\sigma}^2 \right) } \\
\mbox{subject to  } & \lambda\geq 1.
\end{array}\right. \label{MarkovizDROWInternoReform3}
\end{align}
This previous problem can be simplified by analyzing the objective function with respect to $\lambda$. For a fixed $ \beta \in \mathbb{R} $, first note that the function goes to infinity when $\lambda\rightarrow 1^{+}$  or  $\lambda\rightarrow\infty$. Now, its second derivative is given by
$$
\frac{\beta^{2}+  4\beta(\eta-\bar{\zeta}) +4\widehat{\sigma}^2}{2(\lambda-1)^{3}}.
$$
Since $\lambda\geq 1$, the sign of the last expression is determined by the sign of its numerator, which, in terms of  $\beta $, is a polynomial with discriminant given by 
$ \left(\eta-\bar{\zeta}\right)^{2}-\widehat{\sigma}^2$. As a consequence of Cauchy-Schwartz inequality, this discriminant is always negative which implies that the polynomial is always positive. Therefore, the objective function in \eqref{MarkovizDROWInternoReform3} is convex and has a unique minimum value in the region $\lambda\geq 1$. This minimum is reached at $\lambda^{*}=1+\frac{1}{\varepsilon}\sqrt{\frac{\beta^{2}}{4}+\beta(\eta-\bar{\zeta}) +\widehat{\sigma}^2} $ and (\ref{MarkovizDROWInternoReform3}) can be rewritten as
\begin{equation} \label{MarkovizDROWInternoReform4}
\inf_{\beta\in\mathbb{R}}\left(\varepsilon^{2}+ \beta(\eta-\bar{\zeta}) +\widehat{\sigma}^2+\varepsilon\sqrt{ \beta^{2}+  4\beta(\eta-\bar{\zeta}) +4\widehat{\sigma}^2} \right).
\end{equation}
Since the objective function is differentiable for all $\beta\in \mathbb{R} $, after some calculations we obtain that the infimum is attained at $\beta^{*}=2(\bar{\zeta}-\eta)+2(\bar{\zeta}-\eta)\sqrt{\dfrac{\left(\widehat{\sigma}^2 -(\bar{\zeta}-\eta)^2\right)}{\varepsilon^2 -(\bar{\zeta}-\eta)^2}}$ and the optimal value of (\ref{VarianConMeanDROW}) is
$$\left(\sqrt{\widehat{\sigma}^2 -\left(\bar{\zeta}-\eta\right)^{2}  }+\sqrt{\varepsilon^{2}- \left(\bar{\zeta}-\eta\right)^{2} } \right)^{2}.$$
 \endproof

\subsection{Proofs of Theorems \ref{Prop:Reformul1MarkovizDROWReformLargaLargaSupportAcot} and \ref{Prop:Reformul1MarkovizDROWReformLargaLarga} and Corollaries} \label{Sec:Appendix:ProofThmMeanVarianze}

Before proceeding with the proof of Theorems \ref{Prop:Reformul1MarkovizDROWReformLargaLargaSupportAcot} and \ref{Prop:Reformul1MarkovizDROWReformLargaLarga}, we need the following lemma which allows us to express the feasible set of problem (\ref{MarkovizRobust}) in terms of finite dimensional variables. 

\begin{lemma} \label{Lemma:RobustVersionOfMean}
Assume the same setting as in Section \ref{Apendice:DRVarianza}.  
\begin{enumerate}
    \item[(a)] If $p=1$ and $\Xi$ is an interval, then 
    $$  \sup_{ \mathbb{Q}\in \mathcal{B}_{\varepsilon}(\widehat{\mathbb{P}}_{N})}\mathbb{E}_{\mathbb{Q}}[\zeta]=\min\left\{\bar{\zeta}+\varepsilon,\sup_{\xi\in\Xi}\xi\right\} \quad \mbox{and}\quad   \inf_{ \mathbb{Q}\in \mathcal{B}_{\varepsilon}(\widehat{\mathbb{P}}_{N})}\mathbb{E}_{\mathbb{Q}}[\zeta]=\max\left\{\bar{\zeta}-\varepsilon,\inf_{\xi\in\Xi}\xi\right\}.
$$
    \item[(b)] If $p\geq 1$, $\Xi$ is an interval, and $\sup_{\zeta\in\Xi}\zeta=\infty$, then
    $$\sup_{ \mathbb{Q}\in \mathcal{B}_{\varepsilon}(\widehat{\mathbb{P}}_{N})}\mathbb{E}_{\mathbb{Q}}[\zeta] =\bar{\zeta}+\varepsilon\left(\frac{1}{p}+\frac{p-1}{p^{1/(p-1)}}\right)
    \quad \mbox{and}\quad  \inf_{ \mathbb{Q}\in \mathcal{B}_{\varepsilon}(\widehat{\mathbb{P}}_{N})}\mathbb{E}_{\mathbb{Q}}[\zeta]= \bar{\zeta}-\varepsilon\left(\frac{1}{p}+\frac{p-1}{p^{1/(p-1)}}\right).
    $$
\end{enumerate}
\end{lemma}
\proof
It suffices to demonstrate the first equality for each case, as the second is readily obtained by setting $\inf_{ \mathbb{Q}\in \mathcal{B}_{\varepsilon}(\widehat{\mathbb{P}}_{N})}\mathbb{E}_{\mathbb{Q}}[\zeta]=-\sup_{ \mathbb{Q}\in \mathcal{B}_{\varepsilon}(\widehat{\mathbb{P}}_{N})}\mathbb{E}_{\mathbb{Q}}[-\zeta]$. For $p$-Wasserstein distances, by Theorem  \ref{Thm:ReformulacionDROWInterno} we have that
\begin{equation}  \label{ReformulationInf}
 \sup_{ \mathbb{Q}\in \mathcal{B}_{\varepsilon}(\widehat{\mathbb{P}}_{N})}\mathbb{E}_{\mathbb{Q}}[\zeta]=\left\{\begin{array}{lll}
{\displaystyle \inf_{\lambda\geq 0 } } & \lambda\varepsilon^{p}+\frac{1}{N}\sum_{i=1}^{N}s_{i} &\\
\mbox{subject to}  & \sup_{\zeta\in \Xi}\left(\zeta -\lambda\left|\zeta-\widehat{\zeta}_{i} \right|^{p} \right)\leq s_{i} & \forall\:i=1,\ldots,N, 
\end{array}\right.
\end{equation}
In case (a), considering $B:=\sup_{\xi\in\Xi}\xi$ and $A:=\sup_{\xi\in\Xi}\xi$, it is observed that if $B=\infty$, then
$$ \sup_{\zeta\in \Xi}\left(\zeta -\lambda\left|\zeta-\widehat{\zeta}_{i} \right| \right)=\sup_{\zeta\in (A,\infty)}\left(\zeta -\lambda\left|\zeta-\widehat{\zeta}_{i} \right| \right)=\left\{\begin{array}{ll}
    \widehat{\zeta}_{i} & \mbox{ if }\lambda\geq 1,  \\
     \infty &  \mbox{ if }\lambda<1.
\end{array}\right.$$
Therefore, 
$$\sup_{ \mathbb{Q}\in \mathcal{B}_{\varepsilon}(\widehat{\mathbb{P}}_{N})}\mathbb{E}_{\mathbb{Q}}[\zeta]=\bar{\zeta}+\varepsilon=\min\{\bar{\zeta}+\varepsilon,B\} $$
In addition, if $B<\infty$, then
\begin{align}
\sup_{\zeta\in \Xi}\left(\zeta -\lambda\left|\zeta-\widehat{\zeta}_{i} \right| \right)=& \sup_{\zeta\in (A,B)}\left(\zeta -\max_{|z_{i}|\leq \lambda} z_{i}\left(\zeta-\widehat{\zeta}_{i} \right) \right) \nonumber\\
=& \min_{|z_{i}|\leq \lambda}  \sup_{\zeta\in (A,B)}\left(\zeta - z_{i}\left(\zeta-\widehat{\zeta}_{i} \right) \right) \label{Paso:MinMax}\\
=& \min_{|z_{i}|\leq \lambda} \max\left\{(\widehat{\zeta}_{i}-A)z_{i}+A,(\widehat{\zeta}_{i}-B)z_{i}+B\right\} \nonumber\\
=& \left\{\begin{array}{ll}\widehat{\zeta}_{i} & \mbox{ if }\lambda\geq 1, \\ (\widehat{\zeta}_{i}-B)\lambda +B & \mbox{ if }\lambda <1. \end{array}\right. \nonumber
\end{align}
Equality  (\ref{Paso:MinMax})  is guaranteed by  Von Neumann's minimax theorem (see \cite{Bertsekas}). Therefore, we obtain
\begin{align}
\sup_{ \mathbb{Q}\in \mathcal{B}_{\varepsilon}(\widehat{\mathbb{P}}_{N})}\mathbb{E}_{\mathbb{Q}}[\zeta] =& \min\left\{\inf_{\lambda\geq 1}\left(\lambda\varepsilon+\frac{1}{N}\sum_{i=1}^{N}\widehat{\zeta}_{i}\right),\inf_{\lambda< 1}\left(\lambda\varepsilon+\frac{1}{N}\sum_{i=1}^{N}((\widehat{\zeta}_{i}-B)\lambda+B)\right)\right\} \nonumber\\
=& \min\left\{\bar{\zeta}+\varepsilon,\min\{\bar{\zeta}+\varepsilon,B\}\right\} .\nonumber\\
=& \min\left\{\bar{\zeta}+\varepsilon,B\right\} .\nonumber
\end{align}
In case (b), for $p\geq 2$ we have that
\begin{equation*}
    \sup_{\zeta\in \Xi}\left(\zeta -\lambda\left|\zeta-\widehat{\zeta}_{i} \right|^{p} \right)=\widehat{\zeta}_{i}+\frac{p-1}{\lambda^{1/(p-1)}p^{p/(p-1)}}.
\end{equation*}
Therefore, 
\begin{equation*}
    \sup_{ \mathbb{Q}\in \mathcal{B}_{\varepsilon}(\widehat{\mathbb{P}}_{N})}\mathbb{E}_{\mathbb{Q}}[\zeta] =  \inf_{\lambda\geq 0}\left(\lambda\varepsilon+\frac{p-1}{\lambda^{1/(p-1)}p^{p/(p-1)}}+\bar{\zeta}\right)=\bar{\zeta}+\varepsilon\left(\frac{1}{p}+\frac{p-1}{p^{1/(p-1)}}\right).
\end{equation*}
 \qed
\endproof

\proof[Theorem \ref{Prop:Reformul1MarkovizDROWReformLargaLargaSupportAcot}]
The presented outcome follows straightforwardly from Lemma \ref{Lemma:RobustVersionOfMean}. Specifically, Lemma \ref{Lemma:RobustVersionOfMean}-(a) implies that for $p=1$ and $F(x,\Xi)$ as an interval, the optimization problem
\begin{align*}
    \widehat{J}_{N,p,q}^{\mathrm{A}}(\varepsilon) &=\underset{ x\in\mathbb{X}  }{\min} \min\left\{\frac{1}{N}{\displaystyle\sum_{i=1}^{N}}F\left(x,\widehat{\xi}_{i}\right)+\varepsilon\gamma_{x,F},B_{F}(x)\right\},
\end{align*}
is equivalent to (\ref{MarkovizDROWReformLargaSupportAcot}). The proof of (\ref{MarkovizDROWReformLargaSupportNoAcot}), for $p\geq 1$, $F(x,\Xi)$ as an interval, and $\sup_{\xi\in\Xi}F(x,\xi)=\infty$ for each $x\in\mathcal{X}$, can be established in a similar way using Lemma \ref{Lemma:RobustVersionOfMean}-(b).
\qed
\endproof

\proof[Theorem \ref{Prop:Reformul1MarkovizDROWReformLargaLarga}]
By Proposition \ref{Prop:EpsilonFactible} (\ref{MarkovizRobust}) is equivalent to 
\begin{align}
\widehat{J}_{N,p,q}^{A} &:=\underset{ x\in\mathbb{X}  }{\min}  \sup_{\mathbb{Q}\in\mathcal{B}_{\varepsilon \gamma_{x,F}}(\widehat{\mathbb{P}}_{N}^{x}) }  \mathrm{Var}_{\mathbb{Q}}\left[\zeta\right]  \nonumber \\
&= \underset{ x\in\mathbb{X}  }{\min} \sup_{\left(\eta -\frac{1}{N}\sum_{i=1}^{N}\widehat{\zeta}_{i}^{x,F} \right)^{2} \leq \varepsilon^{2} \gamma_{x,F}^{2}} \left\{ \begin{array}{ll} \sup_{\mathbb{Q}\in\mathcal{B}_{\varepsilon  \gamma_{x,F}}(\widehat{\mathbb{P}}_{N}^{x}) } &  \mathrm{Var}_{\mathbb{Q}}\left[\zeta\right] \\ \mbox{subject to} &\mathbb{E}_{\mathbb{Q}}\left[\zeta\right]=\eta. \end{array} \right.  \label{MarkovizDROW} 
\end{align}  
By Theorem \ref{Thm:EsimaVarConMediaConoc}, the maximization problem in (\ref{MarkovizDROW}) which has a variance as its objective function, can be rewritten. Then (\ref{MarkovizDROW}) is equivalent to
\begin{align}
\widehat{J}_{N,p,q}^{A} &=\underset{ x\in\mathbb{X}  }{\min}  \resizebox{.75\hsize}{!}{$\left\{\begin{array}{ll} {\displaystyle\sup_{\eta} }& \left( \sqrt{\frac{1}{N}{\displaystyle\sum_{i=1}^{N}}F\left(x,\widehat{\xi}_{i} \right) ^{2}-\frac{1}{N^{2}}\left({\displaystyle\sum_{i=1}^{N}}F\left( x,\widehat{\xi}_{i}\right)\right)^{2} }+\sqrt{\varepsilon^{2}\gamma_{x,F}^{2}-\left(\eta-\frac{1}{N}{\displaystyle\sum_{i=1}^{N}}F\left( x,\widehat{\xi}_{i}\right)\right)^{2} }  \right)^{2} \\ \mbox{subject to  } &  {\displaystyle\left(\eta -\frac{1}{N}\sum_{i=1}^{N}\widehat{\zeta}_{i}^{x,F} \right)^{2} \leq \varepsilon^{2}\gamma_{x,F}^{2} } \end{array}\right.$} \label{MarkovizDROWReformLarga1}
\end{align}
But, note that the internal maximization problem of (\ref{MarkovizDROWReformLarga1}) can be explicitly solved. Actually, this problem reaches its optimal value in $\eta^{*}=\frac{1}{N}\sum_{i=1}^{N}\widehat{\zeta}_{i}^{x,F}$. Therefore, (\ref{MarkovizDROWReformLarga1}) can be rewritten as
\begin{align}
\widehat{J}_{N,p,q}^{A} (\varepsilon) &=\underset{ x\in\mathbb{X}  }{\mathrm{minimize}} \left( \sqrt{\frac{1}{N}{\displaystyle\sum_{i=1}^{N}}F\left(x,\widehat{\xi}_{i}\right) ^{2}-\frac{1}{N^{2}}\left({\displaystyle\sum_{i=1}^{N}}F\left( x,\widehat{\xi}_{i}\right)\right)^{2} }+\varepsilon\gamma_{x,F} \right)^{2}. \nonumber 
\end{align} 
\qed
\endproof

%*************************************************************************
\subsection{Proof of Proposition \ref{Prop:IgualdadEnReformulaciones}}
For the proof of this proposition, it is necessary to establish the following lemmas.
\begin{lemma} \label{LemmaLipschitz}
    Let $f:\mathbb{R}^{m}\rightarrow\mathbb{R}$ be a convex $2$-Lipschitz function, and $\Theta_{f}:=\left\{z\in\mathbb{R}^{m}\: :\: f^{*}(z)<\infty\right\}$ where $f^{*}$ is the conjugate function of $f$. Then  $\sup_{z\in \Theta_{f}}\|z\|_{2}=\|f\|_{\mathrm{Lip},2}$.
\end{lemma}

\proof
We initiate this demonstration by proving that $\sup_{z\in \Theta_{f}}\|z\|_{2}\leq\|f\|_{\mathrm{Lip},2}$. Let $z\in \Theta_{f}$, which implies that $f^{*}(z)<\infty$. Therefore, there exists a sequence $\left\{z_{n}\right\}_{n=1}^{\infty}$ such that $f^{*}(z)\leq \langle z,z_{n} \rangle -f(z_n)+\frac{1}{n}$ for all $n\geq 1$. This implies that $\langle z,y \rangle -f(y)\leq \langle z,z_{n} \rangle -f(z_{n})+\frac{1}{n}$  for all $y$ and $n\geq 1$, which is equivalent to $0\leq \langle z,z_{n}-y \rangle -(f(z_{n})-f(y))+\frac{1}{n}$ for all $y$ and $n\geq 1$. In particular, making $y=z_{n}+z$ we have $0\leq \langle z,-z \rangle -(f(z_{n})-f(z_{n}+z))+\frac{1}{n}$ for all $n\geq 1$. Therefore, $\|z\|^{2}_{2}\leq -(f(z_{n})-f(z_{n}+z))+\frac{1}{n}$. This implies that
$$\|z\|_{2}\leq -\frac{(f(z_{n})-f(z_{n}+z))}{\|z\|_{2}}+\frac{1}{n\|z\|_{2}}\leq \|f\|_{\mathrm{Lip},2}+\frac{1}{n\|z\|_{2}}.$$
Taking the supremum with respect to $z\in \Theta_{f}$ and the limit when $n$ tends to infinity we obtain $\sup_{z\in \Theta_{f}}\|z\|_{q}\leq\|f\|_{\mathrm{Lip},q}$.

Now we are going to prove the above inequality in the opposite direction. Let $z\in \Theta_{f}$, the we have $\langle z,x \rangle -f(x)\leq f^{*}(z)$ for all $x$. In particular, for all $x$ and $y$ we have
\begin{equation}\label{Eq:LemmaPropLipstInequality1}
f(y)-f(x)\leq \langle z,y-x \rangle+ f(y)+f^{*}(z)-\langle z,y\rangle.
\end{equation}
From this same inequality, we also have 
\begin{equation}\label{Eq:LemmaPropLipstInequality2}
f(x)-f(y)\leq \langle z,x-y \rangle+ f(x)+f^{*}(z)-\langle z,x\rangle.
\end{equation}
Therefore, from (\ref{Eq:LemmaPropLipstInequality1}) it can be inferred
\begin{align*}
    f(y)-f(x) &\leq \inf_{z\in \Theta_{f} } \left(\langle z,y-x\rangle + f(y)+f^{*}(z)-\langle z,y \rangle\right) \\
    &\leq \sup_{z\in \Theta_{f}} \left(\langle z,y-x\rangle+f(y) \right) +\inf_{z\in \Theta_{f}}\left(f^{*}(z)-\langle z,y\rangle\right)\\
    &= \sup_{z\in \Theta_{f}} \left(\langle z,y-x\rangle \right) +f(y)-f^{**}(y)\\
    &= \sup_{z\in \Theta_{f}} \left(\langle z,y-x\rangle \right)\\
    &\leq  \|y-x\|_{2}\sup_{z\in \Theta_{f}} \|z\|_{2}.
\end{align*}
Analogously, we also have 
$$-(f(y)-f(x))\leq \|y-x\|_{2}\sup_{z\in \Theta_{f}} \|z\|_{2}.$$
Combining the last two inequalities we have $\|f\|_{\mathrm{Lip},2}\leq \sup_{z\in \Theta_{f}} \|z\|_{2}$.
\qed
\endproof

The ideas and strategy of the following Lemma proof are taken from the proof of Lemma 47 in \cite{Shafieezadeh-AbadehKuhn2019}. 

\begin{lemma} \label{LemmaSupremoInterno}
    Let $f:\mathbb{R}^{m}\rightarrow\mathbb{R}$ be a convex $2$-Lipschitz function, and $\lambda \geq 0$. Then
    $$
    \sup_{z\in\mathbb{R}^{m}}\left( f(z)-\lambda\|z-y\|_{2}\right)=\left\{\begin{array}{ll} f(x) &\mbox{ if }\:\|f\|_{\mathrm{Lip},2}\leq \lambda\\ \infty &\mbox{ if }\: \|f\|_{\mathrm{Lip},2}> \lambda, \end{array}\right.
    $$
    for all $y$.
\end{lemma}
\proof
Since $f$ is convex, it follows that $f(z)=f^{**}(z)=\sup_{w\in\Theta_{f}}\left(\langle w,z\rangle-f^{*}(w)\right)$. Therefore,
\begin{align}
    \sup_{z\in\mathbb{R}^{m}}\left( f(z)-\lambda\|z-y\|_{2}\right) &= \sup_{z\in\mathbb{R}^{m}}\left( \sup_{w\in\Theta_{f}}\left(\langle w,z\rangle-f^{*}(w)\right)-\lambda\|z-y\|_{2}\right) \nonumber \\
    &=\sup_{z\in\mathbb{R}^{m}}\left( \sup_{w\in\Theta_{f}}\left(\langle w,z\rangle-f^{*}(w)\right)-\sup_{\|\alpha \|_{2}\leq\lambda}\langle \alpha ,z-y\rangle\right) \nonumber \\
    &=\sup_{z\in\mathbb{R}^{m}}\sup_{w\in\Theta_{f}} \inf_{\|\alpha \|_{2}\leq\lambda}\left( \langle w,z\rangle-f^{*}(w)-\langle \alpha ,z-y\rangle\right) \nonumber \\
    &=\sup_{w\in\Theta_{f}} \inf_{\|\alpha \|_{2}\leq\lambda}\sup_{z\in\mathbb{R}^{m}}\left( \langle w,z\rangle-f^{*}(w)-\langle \alpha ,z-y\rangle\right)  \label{Eqn:UsarBersekas}\\
    &=\sup_{w\in\Theta_{f}} \inf_{\|\alpha \|_{2}\leq\lambda}\sup_{z\in\mathbb{R}^{m}}\left( \langle w-\alpha,z\rangle-f^{*}(w)+\langle \alpha ,y\rangle\right) \nonumber\\
    &=\sup_{w\in\Theta_{f}} \inf_{\|\alpha \|_{2}\leq\lambda}\left\{\begin{array}{ll} \langle \alpha,y\rangle-f^{*}(w) &\mbox{ if }\: w=\alpha,\\ \infty &\mbox{ if }w\neq \alpha. \end{array} \right. \nonumber\\
    &=\sup_{w\in\Theta_{f}}\left\{\begin{array}{ll} \langle w,y\rangle-f^{*}(w) &\mbox{ if }\: \|w\|_{2}\leq \lambda,\\ \infty &\mbox{ if }\:\|w\|_{2}> \lambda. \end{array} \right. \nonumber\\
    &=\left\{\begin{array}{ll} \sup_{w\in\Theta_{f}}\left(\langle w,y\rangle-f^{*}(w)\right) &\mbox{ if }\: \sup_{w\in\Theta_{f}}\|w\|_{2}\leq \lambda,\\ \infty &\mbox{ if }\:\sup_{w\in\Theta_{f}}\|w\|_{2}> \lambda. \end{array} \right. \nonumber\\
    &=\left\{\begin{array}{ll} f(y) &\mbox{ if }\: \|f\|_{\mathrm{Lip},2}\leq \lambda,\\ \infty &\mbox{ if }\:\|f\|_{\mathrm{Lip},2}> \lambda. \end{array} \right. \label{Eqn:UsarLemaLipschitz}
\end{align}
Equality (\ref{Eqn:UsarBersekas}) is due to Proposition 5.5.4 in \cite{Bertsekas2009}, and equality (\ref{Eqn:UsarLemaLipschitz}) is by Lemma \ref{LemmaLipschitz} and the fact that $f$ is convex.
\qed
\endproof

This allows us to proceed with the proof of Proposition \ref{Prop:IgualdadEnReformulaciones}.

\proof[ Proposition \ref{Prop:IgualdadEnReformulaciones}]
Let $x\in\mathcal{X}$ be fixed. Note that, by Theorem \ref{Thm:ReformulacionDROWInterno} we have that 
$$\sup_{\mathbb{Q}\in\mathcal{B}_{\varepsilon}\left(\widehat{\mathbb{P}}_{N}\right) }\mathbb{E}_{\mathbb{Q}}[F(x,\xi)] = \inf_{\lambda\geq 0}\left(\lambda \varepsilon^{p} +\frac{1}{N}\sum_{i=1}^{N}\sup_{\xi\in\mathbb{R}^{m}}\left(F(x,\xi)-\lambda \left\|\xi-\widehat{\xi}_{i}\right\|_{q}^{p}  \right)\right).$$
However, for case (i), from Lemma \ref{LemmaSupremoInterno} we know that
$$
\sup_{\xi\in\mathbb{R}^{m}}\left(F(x,\xi)-\lambda \left\|\xi-\widehat{\xi}_{i}\right\|_{2}  \right)=\left\{\begin{array}{ll} F(x,\widehat{\xi}_{i}) &\mbox{ if }\: \gamma_{x,F,2}\leq \lambda,\\ \infty &\mbox{ if }\:\gamma_{x,F,2}> \lambda. \end{array} \right.
$$
Therefore, for case (i), we have
\begin{align}
    \sup_{\mathbb{Q}\in\mathcal{B}_{\varepsilon}\left(\widehat{\mathbb{P}}_{N}\right) }\mathbb{E}_{\mathbb{Q}}[F(x,\xi)] &= \inf_{\lambda\geq \gamma_{x,F,2} }\left(\lambda \varepsilon +\frac{1}{N}\sum_{i=1}^{N}F(x,\widehat{\xi}_{i})\right) \nonumber\\
    &= \varepsilon \gamma_{x,F,2} + \frac{1}{N}\sum_{i=1}^{N}F(x,\widehat{\xi}_{i})\nonumber\\
    &=\sup_{\mathbb{Q}\in\mathcal{B}_{\varepsilon\gamma_{x,F,q}}\left(\widehat{\mathbb{P}}_{N}^{x,F}\right) } \mathbb{E}_{\zeta\sim\mathbb{Q}}\left[\zeta\right]. \label{Eqn:IgualdadEnfoquesCasop1}
\end{align}
Equality (\ref{Eqn:IgualdadEnfoquesCasop1}) is due to Theorem 
\ref{Prop:Reformul1MarkovizDROWReformLargaLargaSupportAcot}(a). Similarly, for case (ii), we have that 
\begin{align}
    \sup_{\xi\in\mathbb{R}^{m}}\left(F(x,\xi)-\lambda \left\|\xi-\widehat{\xi}_{i}\right\|_{2}^{2}  \right)   &= \sup_{\xi\in\mathbb{R}^{m}}\left( I_{x,F}(\xi)-\lambda \left\|\xi-\widehat{\xi}_{i}\right\|_{2}^{2}  \right) \nonumber\\
    &= \sup_{\xi\in\mathbb{R}^{m}}\left( \sup_{z\in \Uptheta_{x,F}}(\langle z,\xi\rangle -I_{x,F}^{*}(z))-\lambda \left\|\xi-\widehat{\xi}_{i}\right\|_{2}^{2}  \right) \label{eqn1:PruebaPropIgualdadEnfoqueP2}\\
    &=\sup_{z\in \Uptheta_{x,F}}\sup_{\xi\in\mathbb{R}^{m}}\left( \langle z,\xi\rangle -I_{x,F}^{*}(z)-\lambda \left\|\xi-\widehat{\xi}_{i}\right\|_{2}^{2}  \right) \nonumber\\
    &=\sup_{z\in \Uptheta_{x,F} } \left( \frac{\|z\|^{2}}{4\lambda }+\langle z, \widehat{\xi}_{i} \rangle-I_{x,F}^{*}(z) \right)  \nonumber \\
    &= \frac{\gamma_{x,F,2}^{2} }{4\lambda}+I_{x,F}(\widehat{\xi}_{i}). \label{eqn1:IgualdadEnfoqueP2}
\end{align}
Equality (\ref{eqn1:PruebaPropIgualdadEnfoqueP2}) is a consequence of the fact that $F(x,\cdot)$ is convex, and equality (\ref{eqn1:IgualdadEnfoqueP2}) is due to the conditions of the case (ii). From this we obtain the following:
\begin{align}
    \sup_{\mathbb{Q}\in\mathcal{B}_{\varepsilon}\left(\widehat{\mathbb{P}}_{N}\right) }\mathbb{E}_{\mathbb{Q}}[F(x,\xi)] &= \inf_{\lambda\geq 0 }\left(\lambda \varepsilon^{2} +\frac{1}{N}\sum_{i=1}^{N}\left(\frac{\gamma_{x,F,2}^{2} }{4\lambda}+I_{x,F}(\widehat{\xi}_{i})\right)\right) \nonumber\\
    & =\inf_{\lambda\geq 0 }\left(\lambda \varepsilon^{2}+\frac{\gamma_{x,F,2}^{2} }{4\lambda} +\frac{1}{N}\sum_{i=1}^{N}I_{x,F}(\widehat{\xi}_{i})\right) \nonumber\\
    &= \varepsilon \gamma_{x,F,2} + \frac{1}{N}\sum_{i=1}^{N}F(x,\widehat{\xi}_{i})\nonumber\\
    &=\sup_{\mathbb{Q}\in\mathcal{B}_{\varepsilon\gamma_{x,F,q}}\left(\widehat{\mathbb{P}}_{N}^{x,F}\right) } \mathbb{E}_{\zeta\sim\mathbb{Q}}\left[\zeta\right]. \label{Eqn:IgualdadEnfoquesCasop2}
\end{align}
Equality (\ref{Eqn:IgualdadEnfoquesCasop2}) is due to Theorem 
\ref{Prop:Reformul1MarkovizDROWReformLargaLargaSupportAcot}(b).
\qed
\endproof

\end{document}